\documentclass[12pt]{article}

\makeatletter
\renewcommand{\@oddfoot}{\hfill \thepage}
\makeatother

\usepackage[T2A]{fontenc}
\usepackage[cp1251]{inputenc}
\usepackage[russian, english]{babel}
\usepackage{mathrsfs}
\usepackage{amssymb}
\usepackage{amsmath,amsthm,amsfonts}
\usepackage[dvipdf]{graphicx}
\usepackage{fancyhdr}

\begin{document}

\begin{center}
{\bf PROBABILITY LAW FOR THE EUCLIDEAN DISTANCE \\ 
     BETWEEN TWO PLANAR RANDOM FLIGHTS} 
\end{center}

\begin{center}
Alexander D. KOLESNIK\\
Institute of Mathematics and Computer Science\\
Academy Street 5, Kishinev 2028, Moldova\\
E-Mail: kolesnik@math.md 
\end{center}

\vskip 0.2cm

\begin{abstract}
We consider two independent symmetric Markov random flights $\bold Z_1(t)$ and $\bold Z_2(t)$ performed by the particles that simultaneously start from the origin of
the Euclidean plane $\Bbb R^2$ in random directions distributed uniformly on the unit circumference $S_1$ and move with constant finite velocities $c_1>0, \; c_2>0$, respectively. The new random directions are taking uniformly on $S_1$ at random time instants that form independent homogeneous Poisson flows of rates $\lambda_1>0, \; \lambda_2>0$. 

The probability distribution function of the Euclidean distance 
$$\rho(t)=\Vert \bold Z_1(t) - \bold Z_2(t) \Vert, \qquad t>0,$$ 
between $\bold Z_1(t)$ and $\bold Z_2(t)$ at arbitrary time instant $t>0$, is obtained.
\end{abstract}

\vskip 0.1cm


\vskip 0.2cm


\section{Introduction}

\numberwithin{equation}{section}

Random flight in the Euclidean space $\Bbb R^m, \; m\ge 2,$ is performed by the stochastic motion of a particle that moves with some finite speed and 
changes, at random time instants, the direction of motion by choosing it on the unit $(m-1)$-dimensional sphere according to some probability 
distribution. Such highly rich stochastic model can generate a lot of particular random walks that might be distinguished following their  
main features:   

- by the velocity (i.e. the speed of motion is constant, or it is a determenistic function depending on space and time variables, 
  or it is a random variable with given distribution);

- by the stochastic flow of the random time instants in which the particle changes its direction (in other words, by the distribution of the time interval $\tau$ 
  between two successive random instants of the flow); 

- by the probability law of choosing the initial and all next random directions; 

- by the dimension of the phase space $\Bbb R^m, \; m\ge 2$. 

\vskip 0.2cm

In the one-dimensional case the particle can take two possible directions (positive and negative) and this is the only difference from the multidimensional model in which  a continuum of directions is assumed. Note that random flights can be treated in a more general context of random evolutions (see, for instance, \cite{pin} and the bibliography therein). 

Random flights are of a special interest due to both their great theoretical importance and, especially, numerous fruitful applications in physics, biology, transport phenomena,  financial modelling and other fields of science and technology. That is why during last decades such stochastic processes have become the subject of extensive researches and a great deal of relevant works were published. 

The most studied one-dimensional stochastic motion at finite speed is represented by the classical Goldstein-Kac telegraph process and its numerous generalizations. At present, this subject numbers several dozens of works. 

As far as the multidimensional models are concerned, their properties are studied to a much less extent. Markov random flight in the Euclidean plane $\Bbb R^2$ with unit  speed and the uniform choice of directions has first been examined by Stadje \cite{sta3}. In this work the explicit form of the transition density of this process 
was derived. The same result (for arbitrary constant speed) was re-obtained then by different methods in \cite{mas}, \cite{kol6}, \cite{kol4}, \cite{kol3}. The similar results for Markov random flights in the Euclidean spaces $\Bbb R^4$ and $\Bbb R^6$ (also with constant speed and the uniform choice of directions) were obtained in \cite{kol5}, \cite{kol3}, \cite{kol2}. The remarkable peculiarity is that in all of these even-dimensional spaces the transition densities are obtained in explicit forms. Moreover, in the spaces $\Bbb R^2$ and $\Bbb R^4$ these densities are surprisingly expressed in terms of elementary function. In the odd-dimensional spaces the analysis is much more difficult. The only result by Stadje \cite{sta2} related to the Markov random flight in the Euclidean space $\Bbb R^3$ with unit speed and the uniform choice of directions yields the transition density in the form of a fairly complicated integral of inverse hyperbolic functions which cannot apparently be evaluated explicitly. 
In all of these works it is supposed that the motions are driven by a homogeneous Poisson process of constant rate. This means that the random time instants, in which the particle changes its direction, form a homogeneous Poisson flow. From this fact it follows that the time interval $\tau$ between two successive turns is an exponentially  distributed random variable. This explains the term "Markov random flight". 

A natural generalization of the Markov random flight is referred to the case when $\tau$ has some other distribution different from the exponential one. Recently a series of works has appeared dealing with the motions when $\tau$ is an Erlang, Dirichlet or Pearson-Dirichlet distributed random variable. In this case the particle changes its direction only at those random time instants of the governing stochastic flow, whose numbers are multiple to some given integer $k\ge 1$, while instants with other numbers are ignored.  Clearly, for $k=1$ we have the Markovian case and for $k=2$ the time interval $\tau$ between turns is an Erlang-distributed random variable.

Le Ca\"er \cite{lecaer1,lecaer2} has examined random broken line in the Euclidean space $\Bbb R^m$ that begins at the origin and consists of uniformly oriented random segments of Dirichlet-distributed random lengths. Such extremely interesting object of stochastic geometry is quite similar to the sample path of Dirichlet random flight. A random broken line can be determined by the ordered triple $(m,n,q)$, where $m$ is the dimension of the space, $n$ stands for the number of break points and $q$ is a positive number  characterizing  the intensity of such points. It was shown in \cite{lecaer1,lecaer2} that for some particular values of parameter $q$ (namely, for  $q=m/2-1, \; m\ge 3$, and $q=m-1, \; q=m, \; m\ge 2,$) the distribution of the end-point of the random broken line can be obtained in an explicit form. Another highly interesting result yields the pairs $(m,n)$ of the dimension $m$ and of the number of break points $n$ under which the end-point has the uniform distribution. This mysterious and hardly explicable effect of appearing the uniform distribution after some changes of direction, first discovered for symmetric Markov random flights in the spaces $\Bbb R^2$ and  $\Bbb R^4$ (see \cite{kol5, kol6}), is turned out to keep its validity also for Dirichlet-distributed uniformy oriented random broken lines. Some of these results were re-obtained in the recent work by Letac and Piccioni \cite{letac} by means of an alternative method based on the Stieltjes transforms of the stochastic broken lines instead of Fourier transforms used in \cite{lecaer1, lecaer2}. 

The integral relations, structurally similar to those obtained in the Markovian case \cite{kol3}, were derived in \cite{pogor2} for Erlang random flight with constant speed and the uniform choice of directions. A symmetric random flight with random velocity and general distribution of $\tau$ was examined in \cite{pogor1} and a renewal equation for the characteristic function of the transition density of the process was given. Some limit theorems for a directionally reinforced random walk were established in  \cite{ghosh}. 

Among a great deal of problems related to the properties of random flights, the problem of finding the probability law for the Euclidean distance between two random flights is of a special interest. This is determined by the importance of such characteristics from the point of view of describing the various kinds of interactions between two moving particles. Such random motions with interaction can serve as very good and adequate mathematical models for describing various
real phenomena in physics, chemistry, biology, financial markets and other fields. For example, in physics and chemistry the
particles can be treated as the atoms or molecules of the substance and their interaction can provoke a physical or chemical reaction. In
biology the particles can be imagined as the biological objects (cells, bacteria, animals etc.) and their "interaction" can mean
creating a new cell (or, contrary, killing the cell), launching an infection mechanism or founding a new animal population, respectively. In
financial markets the moving particles can be interpreted as oscillating exchange rates or stock prices and their "interaction" can mean gaining or ruining.

However, despite the importance of this problem, it has not almost been studied in the literature. The only recent result \cite{kol1} yields a closed-form formula for the probability distribution function of the Euclidean distance between two independent Goldstein-Kac telegraph processes. To the best of our knowledge, the multidimensional counterparts of this problem were not examined yet. In the present article we take the first step in studying this problem in the Euclidean plane $\Bbb R^2$.  

Let $\bold Z_1(t)$ and $\bold Z_2(t)$ denote the positions of two particles on $\Bbb R^2$ at arbitrary time instant $t>0$. As is noted above, in
describing the processes of interaction the crucial role belongs to the Euclidean distance between the particles
$$\rho(t) = \Vert\bold Z_1(t) - \bold Z_2(t) \Vert , \qquad t>0.$$ 
It is quite natural to consider that the particles do not "feel" each
other if $\rho(t)$ is large. In other words, the forces acting
between the particles are negligible if the distance $\rho(t)$ is
sufficiently big. However, as soon as the distance between the
particles becomes less than some given $r>0$, the particles can start interacting with some positive probability.
This means that the occurrence of the random event $\{ \rho(t)< r \} $ is the
necessary (but, maybe, not sufficient) condition for launching the
process of interaction at time instant $t>0$. Therefore, the
distribution $\text{Pr} \{ \rho(t)< r \}$ plays a very important role in analyzing such processes and it is the main objective of our reseach. 

The article is organized as follows. In Section 2 we recall some basic properties of the planar Markov random flight with the uniform choice of  
directions. Three auxiliary lemmas are formulated and proved in Section 3 which we will substantially be relying on. In Section 4 
we prove the principal result of the article representing the probability distribution function for the Euclidean distance between two 
independent planar Markov random flights at arbitrary time instant $t>0$. Although this distribution is expessed in terms of fairly complicated 
integrals, it, nevertheless, can numerically be evaluated by means of the standard package of mathematical programs.

\section{Basic Properties of Symmetric Planar Random Flight}

\numberwithin{equation}{section}

Consider the following stochastic model. A particle is  
located at the origin $\bold 0 = (0, 0)$ of the Euclidean plane $\Bbb R^2$. 
At the initial time instant $t=0$ it starts from $\bold 0$ and moves with constant finite speed $c>0$. 
The initial direction of motion is the two-dimensional random
vector with uniform distribution on the unit circumference
$$S_1 = \left\{ \bold z=(x, y)\in \Bbb R^2: \; \Vert\bold z\Vert^2 = x^2+y^2=1 \right\} .$$
The motion is controlled by a homogeneous Poisson process of rate $\lambda>0$ as follows. 
At the moment of the Poisson event occurrence the particle
instantaneously takes on a new random direction with uniform distribution on $S_1$, 
independently of its previous motion and keeps moving at the same speed $c$ until a new Poisson 
event occurs, then it takes a new random direction, and so on.  

Let $\bold Z(t)=(X(t), Y(t))$ denote the particle's position in the plane $\Bbb R^2$ at an arbitrary time instant $t>0$. 
It is clear that, with probability 1, the process $\bold Z(t)$ is located in the disc of radius $ct$ 
centered at the origin $\bold 0$:  
$$\bold B_{ct} = \left\{ \bold z=(x,y)\in \Bbb R^2: \; \Vert\bold z\Vert^2 = x^2+y^2 \le c^2t^2 \right\} .$$

Let $d\bold z = (dx, dy)$ be the infinitesimal element of the disc $\bold B_{ct}$ with the Lebesgue measure
$\mu(d\bold z) = dx \; dy$. The distribution 
\begin{equation}\label{prop1}
\text{Pr} \left\{ \bold Z(t)\in d\bold z \right\}, \quad
d\bold z\in\bold B_{ct}, \;\; t\ge 0, 
\end{equation}
consists of two components. The singular component of distribution (\ref{prop1}) corresponds to the case
when no one Poisson event occurs in the time interval $(0,t)$ and is concentrated on the circumference
$$S_{ct}=\partial\bold B_{ct} =
\left\{ \bold z=(x,y)\in \Bbb R^2: \; x^2+y^2=c^2t^2 \right\} .$$
In this case the particle is located on $S_{ct}$ and the probability of this event is
\begin{equation}\label{prop2}
\text{Pr} \left\{ \bold Z(t)\in S_{ct} \right\} = e^{-\lambda t} , \qquad t>0. 
\end{equation}

If at least one Poisson event occurs, the particle is located
strictly inside the disc $\bold B_{ct}$, and the probability of this event is
\begin{equation}\label{prop3}
\text{Pr} \left\{ \bold Z(t)\in \text{int} \; \bold B_{ct} \right\} = 1-e^{-\lambda t}, \qquad t>0. 
\end{equation}
The part of distribution (\ref{prop1}) corresponding to this case is concentrated in the interior of disc $\bold B_{ct}$ 
$$\text{int} \; \bold B_{ct} = \left\{ \bold z=(x,y)\in \Bbb R^2: \; \Vert\bold z\Vert^2 = x^2+y^2 < c^2t^2 \right\} ,$$
and forms its absolutely continuous component.

The two-dimensional stochastic process $\bold Z(t)=(X(t), Y(t))$ is referred to as the symmetric planar random flight. 
The principal result states that the transition density $p(\bold z,t)$ of $\bold Z(t)$ at an arbitrary time instant $t>0$ 
is given by the formula:
\begin{equation}\label{prop4}
p(\bold z,t) = \frac{e^{-\lambda t}}{2\pi ct} \; \delta(c^2t^2 - \Vert\bold z\Vert^2) + \frac{\lambda}{2\pi c} 
\frac{\exp \left( -\lambda t + \frac{\lambda}{c} \sqrt{c^2t^2-\Vert\bold z\Vert^2} \right)
}{\sqrt{c^2t^2-\Vert\bold z\Vert^2}} \; \Theta(ct-\Vert\bold z\Vert), 
\end{equation}
$$\bold z\in\bold B_{ct}, \quad \Vert\bold z\Vert^2=x^2+y^2, \quad t>0 ,$$
where $\delta(x)$ is the Dirac delta-function and $\Theta(x)$ is the Heaviside step function. Emphasize that the 
term "density" related to (\ref{prop4}) is treated in the sense of generalized functions. The density (\ref{prop4}) was first derived 
by Stadje \cite{sta3} and re-obtained then by different methods by other authors in \cite{mas}, \cite{kol6}, \cite{kol4}. 

The first term on the right-hand side of (\ref{prop4}) represents the density of the singular component of
distribution (\ref{prop1}), while the second term is the density of the absolutely continuous component of distribution (\ref{prop1}).

The density (\ref{prop4}) has an especially interesting form in polar coordinates: 
\begin{equation}\label{prop5}
\tilde p(r,\alpha,t) = \frac{r e^{-\lambda t}}{2\pi ct} \; \delta(c^2t^2 - r^2) + \frac{\lambda r}{2\pi c} 
\frac{\exp \left( -\lambda t + \frac{\lambda}{c} \sqrt{c^2t^2-r^2} \right)}{\sqrt{c^2t^2-r^2}} \; \Theta(ct-r), 
\end{equation}
$$0<r\le ct, \qquad 0\le\alpha<2\pi, \qquad t>0 .$$  
The radial component $R(t)$, representing the Euclidean distance from the origin of the random point $\bold Z(t)$ 
at an arbitrary time instant $t>0$
$$R(t)= \left\Vert\overrightarrow{(\bold 0, \bold Z(t))} \right\Vert = \sqrt{X^2(t) + Y^2(t)}$$
is independent of its angular component $\alpha(t)$ representing the polar angle between the random vector 
$\overrightarrow{(\bold 0, \bold Z(t))}$ and positive half of the $x$-axis (counter clock-wise circuit). 

It is obvious that, with probability 1, $\; 0<R(t)\le ct$ and $0\le\alpha(t)<2\pi$ for any $t>0$. From (\ref{prop5}) it
follows that the density of $R(t)$ is given by the formula 
\begin{equation}\label{prop6}
f(r,t)=f_{R(t)}(r,t) = \frac{r e^{-\lambda t}}{ct} \delta(ct - r) + \frac{\lambda}{c} \; \frac{r}{\sqrt{c^2t^2 - r^2}} 
\exp\left( -\lambda t + \frac{\lambda}{c} \sqrt{c^2t^2 - r^2} \right) \; \Theta (ct-r), 
\end{equation}
$$r\in(0, ct], \qquad t>0,$$
and the probability distribution function of $R(t)$ has the form
\begin{equation}\label{prop7}
F(r,t) = F_{R(t)}(r,t)=\text{Pr} \{ R(t) < r \} = \left\{ \aligned 0, \qquad\qquad\qquad & \text{if} \; r\in (-\infty, \; 0], \\
1 - \exp\left( -\lambda t + \frac{\lambda}{c}
\sqrt{c^2t^2-r^2} \right), \quad & \text{if} \; r\in (0, \; ct], \\
1, \qquad\qquad\qquad & \text{if} \; r\in (ct, \; +\infty) .
\endaligned \right. 
\end{equation}
In (\ref{prop6}) the first term is the singular part of the density $f(r,t)$ concentrated at the single 
singularity point $r=ct$ while the second term 
\begin{equation}\label{pprop6}
f^{ac}(r,t) = \frac{\lambda}{c} \; \frac{r}{\sqrt{c^2t^2 - r^2}} 
\exp\left( -\lambda t + \frac{\lambda}{c} \sqrt{c^2t^2 - r^2} \right) \; \Theta (ct-r), 
\end{equation}
represents the density of the absolutely continuous part of distribution function (\ref{prop7}) concentrated 
in the open interval $(0, ct)$.

The angular component $\alpha(t)$ of the random vector $\overrightarrow{(\bold 0, \bold Z(t))}$ has the uniform 
distribution in the interval $[0, 2\pi)$ with the density 
\begin{equation}\label{prop8}
f_{\alpha(t)}(\gamma,t) = \frac{1}{2\pi} , \qquad \gamma\in [0, 2\pi), \quad t>0,
\end{equation}
which does not depend on time $t$.

The marginal $X(t)$ of the planar random flight $\bold Z(t)$ (that is, the projection of $\bold Z(t)$ onto $x$-axis) represents a one-dimensional stochastic motion with random speed. The telegraph processes with random velocities have already been examined in the literature (see, for instance, \cite{sta1}), however, one can check that the marginal $X(t)$ is not a telegraph process. The same concerns the second marginal $Y(t)$. Using (\ref{prop4}) one can easily show that $X(t)$ has the density 
\begin{equation}\label{prop9}
\aligned 
p(x,t) & = \frac{\partial}{\partial x} \; \text{Pr} \{ X(t)<x \} \\ 
& = \frac{e^{-\lambda t}}{\pi \sqrt{c^2t^2-x^2}} + \frac{\lambda e^{-\lambda t}}{2c} \biggl[ I_0\left( \frac{\lambda}{c}\sqrt{c^2t^2-x^2} \right) + 
\bold L_0\left( \frac{\lambda}{c}\sqrt{c^2t^2-x^2} \right) \biggr] ,
\endaligned 
\end{equation}
$$x\in (-ct, \; ct), \qquad t>0,$$
where $I_0(z)$ and $\bold L_0(z)$ are the modified Bessel and Struve functions of order zero, respectively, with series representations 
$$I_0(z) = \sum_{k=0}^{\infty} \frac{1}{(k!)^2} \left( \frac{z}{2} \right)^{2k}, \qquad \bold L_0(z) = \sum_{k=0}^{\infty} \frac{1}{\Gamma^2\left( k+\frac{3}{2} \right)} \left( \frac{z}{2} \right)^{2k+1} .$$
We see that, in contrast to the density of the one-dimensional Goldstein-Kac telegraph process (see. for instance, \cite[Section 0.4]{pin}), density (\ref{prop9}) of the marginal $X(t)$ is absolutely continuous and does not contain any singular component. However, one can notice that the second term of density (\ref{prop9}) is structurally resembling to the absolutely continuous part of the density of the Goldstein-Kac telegraph process. The only difference is the presence of Struve function in (\ref{prop9}), while for the telegraph process it is replaced by the term $\frac{1}{\lambda} \; \frac{\partial}{\partial t} I_0\left( \frac{\lambda}{c}\sqrt{c^2t^2-x^2} \right)$. 

The spatial symmetry of the planar random flight $\bold Z(t)$ implies that the second marginal $Y(t)$ has the density $p(y,t)=(\partial/\partial y) \; \text{Pr} \{ Y(t)<y \} $ similar to (\ref{prop9}) with the replacement $x\mapsto y$ everywhere in (\ref{prop9}).

\section{Auxiliary Lemmas}

\numberwithin{equation}{section}

Consider two independent symmetric planar random flights $\bold Z_1(t)$ and $\bold Z_2(t)$ developing with constant 
velocities $c_1>0$ and $c_2>0$, respectively. The evolutions of these random flights are driven by two 
independent Poisson processes $N_1(t)$ and $N_2(t)$ of rates $\lambda_1>0$ and $\lambda_2>0$, respectively, as described in Section 2 above.  

Let $\alpha_1(t)$ and $\alpha_2(t)$ denote the polar angles of the random vectors $\overrightarrow{(\bold 0, \bold Z_1(t))}$ 
and $\overrightarrow{(\bold 0, \bold Z_2(t))}$ with $x$-axis, respectively, at instant $t>0$. As is noted above (see (\ref{prop8})), 
both $\alpha_1(t)$ and $\alpha_2(t)$ have the same uniform distribution in the interval $[0, 2\pi)$ not depending on time $t$. 
Obviously, $\alpha_1(t)$ and $\alpha_2(t)$ are independent for any fixed $t>0$. Introduce the following one-dimensional stochastic process
$$\varphi(t) = \vert\alpha_1(t) - \alpha_2(t)\vert, \qquad t>0,$$
representing the module of difference between these polar anles. It is clear that, with probability 1, $\; 0\le \varphi(t) < 2\pi$ 
for arbitrary $t>0$. In the following lemma we present the distribution of $\varphi(t)$.

\bigskip

{\bf Lemma 1.} {\it The probability distribution function of the process $\varphi(t)$
has the form}
\begin{equation}\label{lem1}
\text{Pr} \{ \varphi(t)<z \} = \left\{ \aligned 0, \qquad\qquad & \text{if} \;\; z\in (-\infty, \; 0], \\
\frac{4\pi z - z^2}{4\pi^2}, \qquad & \text{if} \;\; z\in(0, \;2\pi], \\
1, \qquad\qquad & \text{if} \;\; z\in (2\pi, \; +\infty), \endaligned \right. 
\end{equation}
{\it with the density}
\begin{equation}\label{lem2}
p_{\varphi(t)}(z,t) = \frac{\partial}{\partial z} \; \text{Pr}
\{ \varphi(t)<z \} =  \frac{1}{\pi} - \frac{z}{2\pi^2}, \qquad z\in [0, \;2\pi), \quad t>0.
\end{equation}

\bigskip

{\it Proof.} Since, for any fixed $t>0$, the polar angles $\alpha_1(t)$ and $\alpha_2(t)$
are independent and uniformly distributed in the interval $[0, 2\pi)$ random variables, then the statement of the lemma
immediately emerges by evaluating the probability of the inequality $\vert x-y \vert <z$ in the square 
$[0, 2\pi)\times [0, 2\pi)$ by means of some simple geometric reasonings and the
well-known "two friends meeting problem" of elementary probability theory. $\square$

\bigskip

Note that both the probability distribution function (\ref{lem1}) and density (\ref{lem2}) do not depend on time 
variable $t$.

Consider now the following stochastic process:
$$\eta(t) = \cos(\varphi(t)) = \cos(\vert\alpha_1(t)-
\alpha_2(t)\vert), \qquad t>0.$$
Clearly, $-1\le\eta(t)\le 1$ with probability 1 for any $t>0$. 
The next lemma concerns the distribution of the process $\eta(t)$.

\bigskip

{\bf Lemma 2.} {\it The probability distribution function of the process $\eta(t)$ has
the form}
\begin{equation}\label{lem3}
\text{Pr} \{ \eta(t)<z \} = \left\{ \aligned 0, \qquad\qquad & \text{if} \;\; z\in(-\infty, \; -1], \\
1 - \frac{1}{\pi} \arccos{z}, \qquad & \text{if} \;\; z\in (-1, \; 1], \\
1, \qquad\qquad & \text{if} \;\; z\in (1, \; +\infty), \endaligned \right. 
\end{equation}
 {\it with the density}
\begin{equation}\label{lem4}
p_{\eta(t)}(z,t) = \frac{\partial}{\partial z} \; \text{Pr}
\{ \eta(t) < z \} =  \frac{1}{\pi\sqrt{1-z^2}}, \qquad z\in [-1, \; 1], \quad t>0.
\end{equation}

\bigskip

{\it Proof.} Taking into account that $0\le\varphi(t)<2\pi$ and according to (\ref{lem2}), 
we have for arbitrary $z\in (-1, \; 1]$:  
$$\aligned
\text{Pr} \{ \eta(t)<z \} & = \text{Pr} \{ \cos(\varphi(t))<z \} \\
& = \text{Pr} \{ \arccos{z}<\varphi(t)<2\pi-\arccos{z} \} \\
& = \int\limits_{\arccos{z}}^{2\pi-\arccos{z}} \left(
\frac{1}{\pi} - \frac{z}{2\pi^2} \right) \; dz\\
& = \frac{1}{\pi} \left( 2\pi - 2\arccos{z} \right) -
\frac{1}{2\pi^2} \left( \left. \frac{z^2}{2}
\right\vert_{\arccos{z}}^{2\pi - \arccos{z}} \right) \\
& = 1 - \frac{1}{\pi} \arccos{z} .
\endaligned$$
The lemma is proved. $\square$

\bigskip

From (\ref{lem3}) we also obtain the tail of the probability distribution function: 
\begin{equation}\label{llem3}
\text{Pr} \{ \eta(t) > z \} = \left\{ \aligned 1, \qquad\qquad & \text{if} \;\; z\in(-\infty, \; -1], \\
\frac{1}{\pi} \arccos{z}, \qquad & \text{if} \;\; z\in (-1, \; 1], \\
0, \qquad\qquad & \text{if} \;\; z\in (1, \; +\infty). \endaligned \right. 
\end{equation}

Define now the following stochastic process:
\begin{equation}\label{lem5}
\theta(t) = \left\{ \aligned \varphi(t), \qquad\qquad & \text{if}
\;\;\; \varphi(t)\in (0, \; \pi), \\
2\pi-\varphi(t), \qquad & \text{if} \;\; \varphi(t) \in (\pi, \; 2\pi),
\endaligned \right. \qquad t>0. 
\end{equation}
For $\varphi(t)=0$ and $\varphi(t)=\pi$ the process $\theta(t)$ is undefined. 
The process $\theta(t)$ represents the acute angle between the random
vectors $\overrightarrow{(\bold 0, \bold Z_1(t))}$ and
$\overrightarrow{(\bold 0, \bold Z_2(t))}$ at time instant $t>0$.
Obviously, $0<\theta(t)<\pi$ with probability 1 for any $t>0$. Process
(\ref{lem5}) can also be rewritten in the following form:
\begin{equation}\label{lem6}
\theta(t) = \arccos(\cos(\varphi(t))), \qquad t>0. 
\end{equation}
We are interested in the distribution of $\theta(t)$. One can
intuitively expect that $\theta(t)$ has the uniform distribution
in the interval $(0, \pi)$, although this is not an obvious
statement. In the following lemma we give a rigorous proof of this
fact.

\bigskip

{\bf Lemma 3.} {\it The process $\theta(t)$ is distributed
uniformly in the interval $(0, \pi)$, that is,}
\begin{equation}\label{lem7}
\text{Pr} \{ \theta(t)<z \} = \left\{ \aligned 0, \qquad & \text{if} \;\; z\in (-\infty, \;0], \\
\frac{z}{\pi}, \qquad & \text{if} \;\; z\in (0, \; \pi], \\
1, \qquad & \text{if} \;\; z\in (\pi, \; +\infty) \endaligned \right. \qquad t>0,
\end{equation}
 {\it and the density of distribution} (\ref{lem7}) {\it is}
\begin{equation}\label{lem8}
p_{\theta(t)}(z,t) = \frac{\partial}{\partial z} \; \text{Pr}
\{ \theta(t)<z \} =  \frac{1}{\pi}, \qquad z\in (0, \; \pi), \quad t>0. 
\end{equation}

\bigskip

{\it Proof.} If $z\in (0, \; \pi]$, then, according to (\ref{lem6}) and (\ref{lem3}),
we have
$$\aligned 
\text{Pr} \{ \theta(t)<z \} & = \text{Pr} \{ \arccos(\cos(\varphi(t)))<z \} \\
& = \text{Pr} \{ \cos{z} < \cos(\varphi(t)) \le 1 \} \\
& = 1 - \text{Pr} \{ \cos(\varphi(t)) \le \cos{z} \} \\
& = 1 - \left( 1 - \frac{1}{\pi} \arccos(\cos{z}) \right) \\
& = \frac{z}{\pi},
\endaligned $$
proving (\ref{lem7}). $\square$

\section{Main Result}

\numberwithin{equation}{section}

Now we are able to prove our principal result concerning the distribution of the Euclidean distance
\begin{equation}\label{main1}
\rho(t) = \Vert\bold Z_1(t) - \bold Z_2(t) \Vert , \qquad t>0,
\end{equation}
between two independent random flights $\bold Z_1(t)$ and $\bold Z_2(t)$ performed by two particles that simultaneously start from the origin $\bold 0 = (0, 0)$ 
of the plane $\Bbb R^2$ and move with constant velocities $c_1>0, \; c_2>0$, respectively. Their evolutions are driven by two independent homogeneous Poisson 
processes $N_1(t), \; N_2(t)$ of rates $\lambda_1>0, \; \lambda_2>0$, as is described above. For the sake of definiteness, we suppose that $c_1\ge c_2$ 
(otherwise, one can merely change numeration of the processes).

Our goal is to obtain the probability distribution function 
\begin{equation}\label{mmain1}
\Phi(r,t) = \text{Pr} \left\{ \rho(t) < r \right\} , \qquad t>0,
\end{equation}
of the Euclidean distance (\ref{main1}). It is clear that $0<\rho(t)<(c_1+c_2)t$ with probability 1 for any $t>0$, that is, 
the open interval $(0, \; (c_1+c_2)t)$ is the support of the distribution of process $\rho(t)$. Note that, in contrast to the one-dimensional case 
(see \cite{kol1}), the distribution of the Euclidean distance (\ref{main1}) is absolutely continuous in the interval $(0, \; (c_1+c_2)t)$ 
and does not contain any singular component. This means that probability distribution function (\ref{mmain1}) is continuous for $r\in \Bbb R$ and 
does not have any jumps. 

The form of the probability distribution function $\Phi(r,t)$ is somewhat different in the cases $c_1>c_2$ and $c_1=c_2$.
First, we derive a formula for $\Phi(r,t)$ in the more difficult case $c_1>c_2$. The more simple case $c_1=c_2$ 
will separately be examined at the end of this section. 

The method of obtaining a formula for $\Phi(r,t)$ is different from that used in the one-dimensional case (see \cite{kol1}). While in \cite{kol1} 
the method is based on evaluating the probability the particle to be located in a $r$-neighbourhood of the other one, in the multidimensional case 
such approach is impracticable. Instead, in order to derive the probability distribution function $\Phi(r,t)$, we consider the 
random triangle with the vertices $\bold 0, \; \bold Z_1(t), \; \bold Z_2(t)$. Two sides $\overrightarrow{(\bold 0, \bold Z_i(t))}, \; i=1,2,$ 
of this triangle are the random vectors of the lengths $R_i(t)= \left\Vert\overrightarrow{(\bold 0, \bold Z_i(t))} \right\Vert, \; i=1,2, $ with known
densities $f_i(r,t), \; i=1,2,$ respectively, given by (\ref{prop6}) (or distribution functions $F_i(r,t), \; i=1,2,$  
given by (\ref{prop7})). The random acute angle $\theta(t)$ between these vectors has the uniform distribution in the
interval $(0, \pi)$ (see Lemma 3). Therefore, our aim is to find the distribution of the third side $\rho(t)$ of this random triangle. 

Our first result is given by the following theorem. 

\bigskip

{\bf Theorem 1.} {\it Under the condition $c_1>c_2$, the probability distribution function $\Phi(r,t)$ of the Euclidean distance $\rho(t)$ between 
two independent planar random flights $\bold Z_1(t)$ and $\bold Z_2(t)$ has the form:}
\begin{equation}\label{Phi}
\Phi(r,t) = \left\{ 
\aligned 0, \qquad & \text{if} \; r\in (-\infty, \; 0] ,\\
         G(r,t), \qquad & \text{if} \; r\in (0, \; m(t)] ,\\
         H_k(r,t), \qquad & \text{if} \; r\in (m(t), \; M(t)] ,\\
         Q(r,t), \qquad & \text{if} \; r\in (M(t), \; c_1t] ,\\
         U(r,t), \qquad & \text{if} \; r\in (c_1t, \; (c_1+c_2)t] ,\\
         1, \qquad & \text{if} \; r\in ((c_1+c_2)t, \; +\infty) ,
\endaligned \right. \qquad c_1 > c_2, \;\; t>0, \;\; k=1,2,         
\end{equation}
{\it where} 
\begin{equation}\label{mM}
m(t) = \min\{ (c_1-c_2)t, \; c_2t \} , \qquad M(t) = \max\{ (c_1-c_2)t, \; c_2t \} ,
\end{equation}
{\it and functions $G(r,t), \; H_k(r,t), \; Q(r,t), \; U(r,t)$ are given by the formulas:}
\begin{equation}\label{G}
\aligned
G(r,t) & = 1 - \exp\left( -\lambda_1t + \frac{\lambda_1}{c_1} \sqrt{c_1^2t^2 - r^2} \right) \\
& + \frac{\lambda_1 e^{-(\lambda_1+\lambda_2)t}}{\pi c_1} \int\limits_{c_2t-r}^{c_2t+r} \arccos{ \left( \frac{\xi^2 + (c_2t)^2 - r^2}{2c_2t \xi} \right)}  \frac{\xi}{\sqrt{c_1^2t^2 - \xi^2}} \exp\left( \frac{\lambda_1}{c_1} \sqrt{c_1^2t^2 - \xi^2} \right) d\xi \\ 
& - \frac{\lambda_1 e^{-(\lambda_1+\lambda_2)t}}{c_1} \int\limits_0^r \frac{\xi}{\sqrt{c_1^2t^2 - \xi^2}} \exp\left( \frac{\lambda_1}{c_1} \sqrt{c_1^2t^2 - \xi^2} \right) \exp\left( \frac{\lambda_2}{c_2} \sqrt{c_2^2t^2 - (r-\xi)^2} \right) d\xi \\
& + \frac{1}{\pi} \int\limits_0^r d\zeta \biggl\{ f_2^{ac}(\zeta,t) \int\limits_{r-\zeta}^{r+\zeta} \arccos{ \left( \frac{\xi^2 + \zeta^2 - r^2}{2\xi\zeta} \right)} \; f_1^{ac}(\xi,t) \; d\xi \biggr\} \\ 
& + \frac{1}{\pi} \int\limits_r^{c_2t} d\zeta \biggl\{ f_2^{ac}(\zeta,t) \int\limits_{\zeta-r}^{\zeta+r} \arccos{ \left( \frac{\xi^2 + \zeta^2 - r^2}{2\xi\zeta} \right)} \; f_1^{ac}(\xi,t) \; d\xi \biggr\} , 
\endaligned
\end{equation}

\begin{equation}\label{H1}
\aligned
H_1(r,t) & = 1 - e^{-(\lambda_1+\lambda_2)t} \exp\left( \frac{\lambda_1}{c_1} \sqrt{c_1^2t^2 - (r-c_2t)^2} \right) \\
& + \frac{\lambda_1 e^{-(\lambda_1+\lambda_2)t}}{\pi c_1}  \int\limits_{r-c_2t}^{r+c_2t} \arccos{ \left( \frac{\xi^2 + (c_2t)^2 - r^2}{2c_2t \xi} \right)} \; \frac{\xi}{\sqrt{c_1^2t^2 - \xi^2}} \; \exp\left( \frac{\lambda_1}{c_1} \sqrt{c_1^2t^2 - \xi^2} \right) d\xi \\
& - \frac{\lambda_2 e^{-(\lambda_1+\lambda_2)t}}{c_2} \int\limits_0^{c_2t} \frac{\zeta}{\sqrt{c_2^2t^2 - \zeta^2}} 
\exp\left( \frac{\lambda_2}{c_2} \sqrt{c_2^2t^2 - \zeta^2} \right) \exp\left( \frac{\lambda_1}{c_1} \sqrt{c_1^2t^2 - (r-\zeta)^2} \right) d\zeta \\
& + \frac{1}{\pi} \int\limits_0^{c_2t} d\zeta \biggl\{ f_2^{ac}(\zeta,t) \int\limits_{r-\zeta}^{r+\zeta} \arccos{ \left( \frac{\xi^2 + \zeta^2 - r^2}{2\xi\zeta} \right)} \; f_1^{ac}(\xi,t) \; d\xi \biggr\} , 
\endaligned
\end{equation}
$$\text{\it if} \;\; c_1 > 2c_2,$$

\begin{equation}\label{H2}
\aligned
H_2(r,t) & = 1 - \exp\left( -\lambda_1t + \frac{\lambda_1}{c_1} \sqrt{c_1^2t^2 - r^2} \right) + \frac{e^{-(\lambda_1+\lambda_2)t}}{\pi} \arccos\left( \frac{(c_1t)^2+(c_2t)^2 - r^2}{2c_1c_2t^2} \right) \\
& + \frac{\lambda_1 e^{-(\lambda_1+\lambda_2)t}}{\pi c_1}  \int\limits_{c_2t-r}^{c_1t} \arccos{ \left( \frac{\xi^2 + (c_2t)^2 - r^2}{2c_2t \xi} \right)} \; \frac{\xi}{\sqrt{c_1^2t^2 - \xi^2}} \; \exp\left( \frac{\lambda_1}{c_1} \sqrt{c_1^2t^2 - \xi^2} \right) \; d\xi \\
& + \frac{\lambda_2 e^{-(\lambda_1+\lambda_2)t}}{\pi c_2}  \int\limits_{c_1t-r}^{c_2t} \arccos{ \left( \frac{\xi^2 + (c_1t)^2 - r^2}{2c_1t \xi} \right)} \; \frac{\xi}{\sqrt{c_2^2t^2 - \xi^2}} \; \exp\left( \frac{\lambda_2}{c_2} \sqrt{c_2^2t^2 - \xi^2} \right) \; d\xi \\
& - \frac{\lambda_1}{c_1} \; e^{-(\lambda_1+\lambda_2)t} \int\limits_0^r \frac{\xi}{\sqrt{c_1^2t^2 - \xi^2}} \exp\left( \frac{\lambda_1}{c_1} \sqrt{c_1^2t^2 - \xi^2} \right) \exp\left( \frac{\lambda_2}{c_2} \sqrt{c_2^2t^2 - (r-\xi)^2} \right) d\xi \\
& + \frac{1}{\pi} \int\limits_0^r d\zeta \biggl\{ f_2^{ac}(\zeta,t) \int\limits_{r-\zeta}^{r+\zeta} \arccos{ \left( \frac{\xi^2 + \zeta^2 - r^2}{2\xi\zeta} \right)} \; f_1^{ac}(\xi,t) \; d\xi \biggr\} \\ 
& + \frac{1}{\pi} \int\limits_r^{c_1t-r} d\zeta \biggl\{ f_2^{ac}(\zeta,t) \int\limits_{\zeta-r}^{\zeta+r} \arccos{ \left( \frac{\xi^2 + \zeta^2 - r^2}{2\xi\zeta} \right)} \; f_1^{ac}(\xi,t) \; d\xi \biggr\} \\
& + \frac{1}{\pi} \int\limits_{c_1t-r}^{c_2t} d\zeta \biggl\{ f_2^{ac}(\zeta,t) \int\limits_{\zeta-r}^{c_1t} \arccos{ \left( \frac{\xi^2 + \zeta^2 - r^2}{2\xi\zeta} \right)} \; f_1^{ac}(\xi,t) \; d\xi \biggr\} ; 
\endaligned
\end{equation}
$$\text{\it if} \;\; c_1 < 2c_2,$$

\begin{equation}\label{Q}
\aligned
Q(r,t) & = 1 - e^{-(\lambda_1+\lambda_2)t} \left[ \exp\left( \frac{\lambda_1}{c_1} \sqrt{c_1^2t^2 - (r-c_2t)^2} \right) - \frac{1}{\pi} \arccos\left( \frac{(c_1t)^2+(c_2t)^2 - r^2}{2c_1c_2t^2} \right) \right] \\
& + \frac{\lambda_1 e^{-(\lambda_1+\lambda_2)t}}{\pi c_1}  \int\limits_{r-c_2t}^{c_1t} \arccos{ \left( \frac{\xi^2 + (c_2t)^2 - r^2}{2c_2t \xi} \right)}  \frac{\xi}{\sqrt{c_1^2t^2 - \xi^2}} \exp\left( \frac{\lambda_1}{c_1} \sqrt{c_1^2t^2 - \xi^2} \right) d\xi \\
& + \frac{\lambda_2 e^{-(\lambda_1+\lambda_2)t}}{\pi c_2}  \int\limits_{c_1t-r}^{c_2t} \arccos{ \left( \frac{\xi^2 + (c_1t)^2 - r^2}{2c_1t \xi} \right)}  \frac{\xi}{\sqrt{c_2^2t^2 - \xi^2}} \exp\left( \frac{\lambda_2}{c_2} \sqrt{c_2^2t^2 - \xi^2} \right) d\xi \\
& - \frac{\lambda_2}{c_2} e^{-(\lambda_1+\lambda_2)t} \int\limits_0^{c_2t} \frac{\zeta}{\sqrt{c_2^2t^2 - \zeta^2}} 
\exp\left( \frac{\lambda_2}{c_2} \sqrt{c_2^2t^2 - \zeta^2} \right) \exp\left( \frac{\lambda_1}{c_1} \sqrt{c_1^2t^2 - (r-\zeta)^2} \right) d\zeta \\
& + \frac{1}{\pi} \int\limits_0^{c_1t-r} d\zeta \biggl\{ f_2^{ac}(\zeta,t) \int\limits_{r-\zeta}^{r+\zeta} \arccos{ \left( \frac{\xi^2 + \zeta^2 - r^2}{2\xi\zeta} \right)} \; f_1^{ac}(\xi,t) \; d\xi \biggr\} \\
& + \frac{1}{\pi} \int\limits_{c_1t-r}^{c_2t} d\zeta \biggl\{ f_2^{ac}(\zeta,t) \int\limits_{r-\zeta}^{c_1t} \arccos{ \left( \frac{\xi^2 + \zeta^2 - r^2}{2\xi\zeta} \right)} \; f_1^{ac}(\xi,t) \; d\xi \biggr\} , 
\endaligned
\end{equation}

\begin{equation}\label{U}
\aligned
U(r,t) & = 1 - e^{-(\lambda_1+\lambda_2)t} \left[ \exp\left( \frac{\lambda_2}{c_2} \sqrt{c_2^2t^2 - (r-c_1t)^2} \right) - \frac{1}{\pi} \arccos\left( \frac{(c_1t)^2+(c_2t)^2 - r^2}{2c_1c_2t^2} \right) \right] \\
& + \frac{\lambda_1 e^{-(\lambda_1+\lambda_2)t}}{\pi c_1}  \int\limits_{r-c_2t}^{c_1t} \arccos{ \left( \frac{\xi^2 + (c_2t)^2 - r^2}{2c_2t \xi} \right)}  \frac{\xi}{\sqrt{c_1^2t^2 - \xi^2}} \exp\left( \frac{\lambda_1}{c_1} \sqrt{c_1^2t^2 - \xi^2} \right) d\xi \\
& + \frac{\lambda_2 e^{-(\lambda_1+\lambda_2)t}}{\pi c_2}  \int\limits_{r-c_1t}^{c_2t} \arccos{ \left( \frac{\xi^2 + (c_1t)^2 - r^2}{2c_1t \xi} \right)}  \frac{\xi}{\sqrt{c_2^2t^2 - \xi^2}} \exp\left( \frac{\lambda_2}{c_2} \sqrt{c_2^2t^2 - \xi^2} \right) d\xi \\
& - \frac{\lambda_1}{c_1} \; e^{-(\lambda_1+\lambda_2)t} \int\limits_{r-c_2t}^{c_1t} \frac{\xi}{\sqrt{c_1^2t^2 - \xi^2}} 
\exp\left( \frac{\lambda_1}{c_1} \sqrt{c_1^2t^2 - \xi^2} \right) \exp\left( \frac{\lambda_2}{c_2} \sqrt{c_2^2t^2 - (r-\xi)^2} \right) d\xi \\
& + \frac{1}{\pi} \int\limits_{r-c_2t}^{c_1t} d\xi \biggl\{ f_1^{ac}(\xi,t) \int\limits_{r-\xi}^{c_2t} \arccos{ \left( \frac{\xi^2 + \zeta^2 - r^2}{2\xi\zeta} \right)} \; f_2^{ac}(\zeta,t) \; d\zeta \biggr\} ,
\endaligned
\end{equation}
{\it where $f_i^{ac}(z,t), \; i=1,2,$ is the absolutely continuous part of the density of process $R_i(t)$ given by formula} (\ref{pprop6}), {\it that is,}
\begin{equation}\label{ContDens}
f_i^{ac}(z,t) = \frac{\lambda_i}{c_i} \; \frac{z}{\sqrt{c_i^2t^2 - z^2}} 
\exp\left( -\lambda_i t + \frac{\lambda_i}{c_i} \sqrt{c_i^2t^2 - z^2} \right) , \qquad 0<z<c_it, \;\; i=1,2.
\end{equation}

\bigskip  

{\it If $c_1=2c_2$, then the probability distribution function $\Phi(r,t)$ has the form:}
\begin{equation}\label{mainA3}
\Phi(r,t) = \left\{ 
\aligned 0, \qquad & \text{if} \; r\in (-\infty, \; 0] ,\\
         G(r,t), \qquad & \text{if} \; r\in (0, \; c_2t] ,\\
         Q(r,t), \qquad & \text{if} \; r\in (c_2t, \; c_1t] ,\\
         U(r,t), \qquad & \text{if} \; r\in (c_1t, \; (c_1+c_2)t] ,\\
         1, \qquad & \text{if} \; r\in ((c_1+c_2)t, \; +\infty) ,
\endaligned \right.  \qquad c_1=2c_2, \;\; t>0,         
\end{equation}
{\it where functions $G(r,t), \; Q(r,t), \; U(r,t)$ are given by formulas} (\ref{G}), (\ref{Q}), (\ref{U}), {\it respectively}.

\bigskip

{\it Proof.} First of all, we note that since $0<\rho(t)<(c_1+c_2)t$ with probability 1 for any $t>0$, then 
\begin{equation}\label{main3}
\aligned 
& \text{Pr} \left\{ \rho(t) < r \right\} = 0, \qquad \text{if} \;\; r\in (-\infty, \; 0] , \\  
& \text{Pr} \left\{ \rho(t) < r \right\} = 1, \qquad \text{if} \;\; r\in ((c_1+c_2)t, \; +\infty) .
\endaligned   
\end{equation}

Let now $r\in (0, \; (c_1+c_2)t]$. Passing to joint distributions, we can write down: 
\begin{equation}\label{main4}
\aligned
\Phi(r,t) & = \text{Pr} \left\{ \rho(t) < r, \; N_1(t)=0, \; N_2(t)=0  \right\} \\
& \; + \text{Pr} \left\{ \rho(t) < r, \; N_1(t)\ge 1, \; N_2(t)=0 \right\} \\
& \; + \text{Pr} \left\{ \rho(t) < r, \; N_1(t)=0, \; N_2(t)\ge 1 \right\} \\
& \; + \text{Pr} \left\{ \rho(t) < r, \; N_1(t)\ge 1, \; N_2(t)\ge 1 \right\} .
\endaligned 
\end{equation}
Let us evaluate separately joint probabilities on the right-hand side of (\ref{main4}). 

\vskip 0.2cm

$\bullet$ {\it Evaluation of} $\; \text{Pr} \left\{ \rho(t) < r, \; N_1(t)=0, \; N_2(t)=0 \right\}$. We note that 
the following equalities for random events hold: 
$$\left\{ N_1(t)=0 \right\} = \left\{ \bold Z_1(t)\in S_{c_1t} \right\} = \left\{ R_1(t)=c_1t \right\}, $$
$$\left\{ N_2(t)=0 \right\} = \left\{ \bold Z_2(t)\in S_{c_2t} \right\} = \left\{ R_2(t)=c_2t \right\} . $$
Then, taking into account that $\cos(\theta(t))=\cos(\varphi(t))=\eta(t)$ (in distribution)
and using (\ref{llem3}), we have for the first joint distribution in (\ref{main4}):
$$\aligned
\text{Pr} & \left\{ \rho(t) < r, \; N_1(t)=0, \; N_2(t)=0 \right\} \\
& = e^{-(\lambda_1+\lambda_2)t} \; \text{Pr} \left\{ R_1^2(t)+R_2^2(t) - 2R_1(t)R_2(t) \cos(\theta(t)) < r^2 \; \bigl\vert \; R_1(t)=c_1t, \; R_2(t)=c_2t \right\} \\
& = e^{-(\lambda_1+\lambda_2)t} \; \text{Pr} \left\{ (c_1t)^2+(c_2t)^2 - 2c_1c_2t^2 \cos(\theta(t)) < r^2 \right\} \\
& = e^{-(\lambda_1+\lambda_2)t} \; \text{Pr} \left\{ \eta(t) > \frac{(c_1t)^2+(c_2t)^2 - r^2}{2c_1c_2t^2} \right\} \\
\\
& = \left\{ \aligned e^{-(\lambda_1+\lambda_2)t}, \qquad\qquad\qquad\qquad  & \text{if} \;\; \frac{(c_1t)^2+(c_2t)^2 - r^2}{2c_1c_2t^2} \le -1,\\
\frac{e^{-(\lambda_1+\lambda_2)t}}{\pi} \arccos\left( \frac{(c_1t)^2+(c_2t)^2 - r^2}{2c_1c_2t^2} \right), \qquad & \text{if} \;\;
-1 < \frac{(c_1t)^2+(c_2t)^2 - r^2}{2c_1c_2t^2} \le 1,\\
0, \qquad\qquad\qquad\qquad & \text{if} \;\; \frac{(c_1t)^2+(c_2t)^2 - r^2}{2c_1c_2t^2} > 1,
\endaligned \right. \endaligned$$

\begin{equation}\label{main5}
= \left\{ \aligned 0, \qquad\qquad\qquad\qquad & \text{if} \;\; r\in (-\infty, \; (c_1-c_2)t],\\
\frac{e^{-(\lambda_1+\lambda_2)t}}{\pi} \arccos\left( \frac{(c_1t)^2+(c_2t)^2 - r^2}{2c_1c_2t^2} \right), \qquad & \text{if} \;\;
r\in ((c_1-c_2)t, \; (c_1+c_2)t], \\
e^{-(\lambda_1+\lambda_2)t}, \qquad\qquad\qquad\qquad & \text{if} \;\; r\in ((c_1+c_2)t, \; +\infty) \endaligned \right.
\end{equation}

Formula (\ref{main5}) yields the first joint distribution in (\ref{main4}) related to the case when no one Poisson event occurs up to time instant $t>0$ 
and, therefore, the randon points $\bold Z_1(t), \; \bold Z_2(t)$ are located on the spheres $S_{c_1t}, \; S_{c_2t}$, respectively.  

\vskip 0.2cm

$\bullet$ {\it Evaluation of} $\; \text{Pr} \left\{ \rho(t) < r, \; N_1(t)\ge 1, \; N_2(t)=0 \right\}$. Note that   
$$\left\{ N_1(t)\ge 1 \right\} = \left\{ \bold Z_1(t)\in \text{int} \; \bold B_{c_1t} \right\} = \left\{ R_1(t)\in (0, \; c_1t) \right\}, $$
$$\left\{ N_2(t)=0 \right\} = \left\{ \bold Z_2(t)\in S_{c_2t} \right\} = \left\{ R_2(t)=c_2t \right\} .$$
Taking into account that $\cos(\theta(t))=\cos(\varphi(t))=\eta(t)$ (in distribution), we have for the second joint distribution in (\ref{main4}):
\begin{equation}\label{main6}
\aligned 
& \text{Pr} \left\{ \rho(t) < r, \; N_1(t)\ge 1, \; N_2(t)=0 \right\} \\
& = e^{-\lambda_2t} \; \text{Pr} \biggl\{ R_1^2(t) + R_2^2(t) - 2R_1(t) R_2(t) \cos(\theta(t)) < r^2, \; R_1(t)\in (0, \; c_1t) \; \bigl\vert \; R_2(t)=c_2t \biggr\} \\
& = e^{-\lambda_2t} \; \text{Pr} \biggl\{ R_1^2(t) + (c_2t)^2 - 2c_2t R_1(t) \cos(\theta(t)) < r^2, \; R_1(t)\in (0, \; c_1t) \biggr\} \\
& = e^{-\lambda_2t} \; \text{Pr} \biggl\{ \eta(t) > \frac{R_1^2(t) + (c_2t)^2 - r^2}{2c_2t R_1(t)} , \; R_1(t)\in (0, \; c_1t) \biggr\} \\
& = e^{-\lambda_2t} \; \int\limits_0^{c_1t} \text{Pr} \biggl\{ \eta(t) > \frac{\xi^2 + (c_2t)^2 - r^2}{2c_2t \xi} \; \biggl\vert \; R_1(t)=\xi \biggr\} \; 
\text{Pr} \bigl\{ R_1(t)\in d\xi \bigr\} . 
\endaligned
\end{equation}
According to (\ref{llem3})
$$\aligned 
\text{Pr} & \biggl\{ \eta(t) > \frac{\xi^2 + (c_2t)^2 - r^2}{2c_2t \xi} \; \biggl\vert \; R_1(t)=\xi \biggr\} \\
& = \left\{ \aligned 1, \qquad\qquad & \text{if} \;\; \frac{\xi^2 + (c_2t)^2 - r^2}{2c_2t \xi} \in(-\infty, \; -1] \; \text{and} \; \xi\in (0, \; c_1t), \\
\frac{1}{\pi} \arccos{ \left( \frac{\xi^2 + (c_2t)^2 - r^2}{2c_2t \xi} \right)}, \quad & \text{if} \;\; \frac{\xi^2 + (c_2t)^2 - r^2}{2c_2t \xi} \in (-1, \; 1] \; \text{and} \; \xi\in (0, \; c_1t), \\
0, \qquad\qquad & \text{if} \;\; \frac{\xi^2 + (c_2t)^2 - r^2}{2c_2t \xi} \in (1, \; +\infty) \; \text{and} \; \xi\in (0, \; c_1t). \endaligned \right.
\endaligned$$
This probability depends on the value of $r$. One can check that, as $\xi\in (0, \; c_1t)$, we get the equalities: 
 
\begin{equation}\label{main7}
\aligned 
& \text{Pr} \biggl\{ \eta(t) > \frac{\xi^2 + (c_2t)^2 - r^2}{2c_2t \xi} \; \biggl\vert \; R_1(t)=\xi \biggr\} \\
\\
& = \left\{ \aligned \frac{1}{\pi} \arccos{ \left( \frac{\xi^2 + (c_2t)^2 - r^2}{2c_2t \xi} \right)}, \qquad & \text{if} \;\; \xi\in (c_2t-r, \; \beta(r)), \\
0, \qquad\qquad & \text{otherwise} , \endaligned \right.  
\endaligned
\end{equation}
$$\text{for} \;\; r\in (0, \; c_2t] ,$$
and 
\begin{equation}\label{main8}
\aligned 
& \text{Pr} \biggl\{ \eta(t) > \frac{\xi^2 + (c_2t)^2 - r^2}{2c_2t \xi} \; \biggl\vert \; R_1(t)=\xi \biggr\} \\
\\
& = \left\{ \aligned 1, \qquad\qquad & \text{if} \;\; \xi\in(0, \; r-c_2t] ,\\
\frac{1}{\pi} \arccos{ \left( \frac{\xi^2 + (c_2t)^2 - r^2}{2c_2t \xi} \right)}, \qquad & \text{if} \;\; \xi\in (r-c_2t, \; \beta(r)), \\
0, \qquad\qquad & \text{otherwise} , \endaligned \right. 
\endaligned
\end{equation}
$$\text{for} \;\; r\in (c_2t, \; (c_1+c_2)t] ,$$
where
$$\beta(r)= \min\bigl\{ c_2t+r, \; c_1t \bigr\}.$$

Substituting (\ref{main7}) into (\ref{main6}) and using (\ref{pprop6}), we obtain for $r\in (0, \; c_2t]$: 
\begin{equation}\label{main9}
\aligned 
& \text{Pr} \left\{ \rho(t) < r, \; N_1(t)\ge 1, \; N_2(t)=0 \right\} \\
& = \frac{\lambda_1 e^{-(\lambda_1+\lambda_2)t}}{\pi c_1}  \int\limits_{c_2t-r}^{\beta(r)} \arccos{ \left( \frac{\xi^2 + (c_2t)^2 - r^2}{2c_2t \xi} \right)} \; \frac{\xi}{\sqrt{c_1^2t^2 - \xi^2}} \; \exp\left( \frac{\lambda_1}{c_1} \sqrt{c_1^2t^2 - \xi^2} \right) \; d\xi ,
\endaligned
\end{equation}
$$\text{for} \;\; r\in (0, \; c_2t] .$$

Similarly, substituting (\ref{main8}) into (\ref{main6}) and using (\ref{pprop6}), we have for $r\in (c_2t, \; (c_1+c_2)t]$: 
$$\aligned 
& \text{Pr} \left\{ \rho(t) < r, \; N_1(t)\ge 1, \; N_2(t)=0 \right\} \\
& = \frac{\lambda_1 e^{-(\lambda_1+\lambda_2)t}}{c_1} \int\limits_0^{r-c_2t} \frac{\xi}{\sqrt{c_1^2t^2 - \xi^2}} \; \exp\left( \frac{\lambda_1}{c_1} \sqrt{c_1^2t^2 - \xi^2} \right) \; d\xi \\
& + \frac{\lambda_1 e^{-(\lambda_1+\lambda_2)t}}{\pi c_1}  \int\limits_{r-c_2t}^{\beta(r)} \arccos{ \left( \frac{\xi^2 + (c_2t)^2 - r^2}{2c_2t \xi} \right)} \; \frac{\xi}{\sqrt{c_1^2t^2 - \xi^2}} \; \exp\left( \frac{\lambda_1}{c_1} \sqrt{c_1^2t^2 - \xi^2} \right) \; d\xi ,
\endaligned$$
$$\text{for} \;\; r\in (c_2t, \; (c_1+c_2)t] .$$
Taking into account that 
\begin{equation}\label{intAB}
\int \frac{x}{\sqrt{p^2-x^2}} \; e^{q\sqrt{p^2-x^2}} \; dx = -\frac{1}{q} \; e^{q\sqrt{p^2-x^2}} , \qquad q\neq 0, \quad |x|\le p,
\end{equation}
we finally obtain: 
\begin{equation}\label{main10}
\aligned 
& \text{Pr} \left\{ \rho(t) < r, \; N_1(t)\ge 1, \; N_2(t)=0 \right\} \\
& = e^{-\lambda_2t} \left[ 1 - \exp\left( -\lambda_1t+\frac{\lambda_1}{c_1} \sqrt{c_1^2t^2 - (r-c_2t)^2} \right) \right] \\
& + \frac{\lambda_1 e^{-(\lambda_1+\lambda_2)t}}{\pi c_1}  \int\limits_{r-c_2t}^{\beta(r)} \arccos{ \left( \frac{\xi^2 + (c_2t)^2 - r^2}{2c_2t \xi} \right)} \; \frac{\xi}{\sqrt{c_1^2t^2 - \xi^2}} \; \exp\left( \frac{\lambda_1}{c_1} \sqrt{c_1^2t^2 - \xi^2} \right) \; d\xi ,
\endaligned
\end{equation}
$$\text{for} \;\; r\in (c_2t, \; (c_1+c_2)t] .$$

\vskip 0.2cm

$\bullet$ {\it Evaluation of} $\; \text{Pr} \left\{ \rho(t) < r, \; N_1(t)=0, \; N_2(t)\ge 1 \right\}$. Since 
$$\left\{ N_1(t)=0 \right\} = \left\{ \bold Z_1(t)\in S_{c_1t} \right\} = \left\{ R_1(t)=c_1t \right\} ,$$    
$$\left\{ N_2(t)\ge 1 \right\} = \left\{ \bold Z_2(t)\in \text{int} \; \bold B_{c_2t} \right\} = \left\{ R_2(t)\in (0, \; c_2t) \right\}, $$
then, similarly as above, we get 
$$\aligned 
\text{Pr} & \left\{ \rho(t) < r, \; N_1(t)=0, \; N_2(t)\ge 1 \right\} \\
& = e^{-\lambda_1t} \int\limits_0^{c_2t} \text{Pr} \biggl\{ \eta(t) > \frac{\xi^2 + (c_1t)^2 - r^2}{2c_1t \xi} \; \biggl\vert \; R_2(t)=\xi \biggr\} \; 
\text{Pr} \bigl\{ R_2(t)\in d\xi \bigr\} .
\endaligned$$
According to (\ref{llem3}), the conditional probability in the integrand is:
$$\aligned 
\text{Pr} & \biggl\{ \eta(t) > \frac{\xi^2 + (c_1t)^2 - r^2}{2c_1t \xi} \; \biggl\vert \; R_2(t)=\xi \biggr\} \\
& = \left\{ \aligned 1, \qquad\qquad & \text{if} \;\; \frac{\xi^2 + (c_1t)^2 - r^2}{2c_1t \xi} \in(-\infty, \; -1] \; \text{and} \; \xi\in (0, \; c_2t), \\
\frac{1}{\pi} \arccos{ \left( \frac{\xi^2 + (c_1t)^2 - r^2}{2c_1t \xi} \right)}, \quad & \text{if} \;\; \frac{\xi^2 + (c_1t)^2 - r^2}{2c_1t \xi} \in (-1, \; 1] \; \text{and} \; \xi\in (0, \; c_2t), \\
0, \qquad\qquad & \text{if} \;\; \frac{\xi^2 + (c_1t)^2 - r^2}{2c_1t \xi} \in (1, \; +\infty) \; \text{and} \; \xi\in (0, \; c_2t). \endaligned \right.
\endaligned$$
This formula splits in the following three cases:
$$\text{Pr} \biggl\{ \eta(t) > \frac{\xi^2 + (c_1t)^2 - r^2}{2c_1t \xi} \; \biggl\vert \; R_2(t)=\xi \biggr\} = 0 , \qquad \text{for} \; r\in (0, \; (c_1-c_2)t],$$
\vskip 0.2cm 
$$\aligned 
& \text{Pr} \biggl\{ \eta(t) > \frac{\xi^2 + (c_1t)^2 - r^2}{2c_1t \xi} \; \biggl\vert \; R_2(t)=\xi \biggr\} \\
\\
& = \left\{ \aligned \frac{1}{\pi} \arccos{ \left( \frac{\xi^2 + (c_1t)^2 - r^2}{2c_1t \xi} \right)}, \qquad & \text{if} \;\; \xi\in (c_1t-r, \; c_2t), \\
0, \qquad\qquad & \text{otherwise} , \endaligned \right.  
\endaligned$$
$$\text{for} \;\; r\in ((c_1-c_2)t, \; c_1t] ,$$
and 
$$\aligned 
& \text{Pr} \biggl\{ \eta(t) > \frac{\xi^2 + (c_1t)^2 - r^2}{2c_1t \xi} \; \biggl\vert \; R_2(t)=\xi \biggr\} \\
\\
& = \left\{ \aligned 1, \qquad\qquad & \text{if} \;\; \xi\in(0, \; r-c_1t] ,\\
\frac{1}{\pi} \arccos{ \left( \frac{\xi^2 + (c_1t)^2 - r^2}{2c_1t \xi} \right)}, \qquad & \text{if} \;\; \xi\in (r-c_1t, \; c_2t], \\
0, \qquad\qquad & \text{otherwise} , \endaligned \right. 
\endaligned$$
$$\text{for} \;\; r\in (c_1t, \; (c_1+c_2)t] .$$

Taking into account (\ref{pprop6}), we therefore obtain: 
\begin{equation}\label{mmain11}
\text{Pr} \left\{ \rho(t) < r, \; N_1(t)=0, \; N_2(t)\ge 1 \right\} = 0, \qquad \text{for} \;\; r\in (0, \; (c_1-c_2)t] ,
\end{equation}
\begin{equation}\label{main11}
\aligned 
& \text{Pr} \left\{ \rho(t) < r, \; N_1(t)=0, \; N_2(t)\ge 1 \right\} \\
& = \frac{\lambda_2 e^{-(\lambda_1+\lambda_2)t}}{\pi c_2}  \int\limits_{c_1t-r}^{c_2t} \arccos{ \left( \frac{\xi^2 + (c_1t)^2 - r^2}{2c_1t \xi} \right)} \; \frac{\xi}{\sqrt{c_2^2t^2 - \xi^2}} \; \exp\left( \frac{\lambda_2}{c_2} \sqrt{c_2^2t^2 - \xi^2} \right) \; d\xi ,
\endaligned
\end{equation}
$$\text{for} \;\; r\in ((c_1-c_2)t, \; c_1t] ,$$
and  
\begin{equation}\label{main12}
\aligned 
& \text{Pr} \left\{ \rho(t) < r, \; N_1(t)=0, \; N_2(t)\ge 1 \right\} \\
& = e^{-\lambda_1t} \left[ 1 - \exp\left( -\lambda_2t+\frac{\lambda_2}{c_2} \sqrt{c_2^2t^2 - (r-c_1t)^2} \right) \right] \\
& + \frac{\lambda_2 e^{-(\lambda_1+\lambda_2)t}}{\pi c_2}  \int\limits_{r-c_1t}^{c_2t} \arccos{ \left( \frac{\xi^2 + (c_1t)^2 - r^2}{2c_1t \xi} \right)} \; \frac{\xi}{\sqrt{c_2^2t^2 - \xi^2}} \; \exp\left( \frac{\lambda_2}{c_2} \sqrt{c_2^2t^2 - \xi^2} \right) \; d\xi ,
\endaligned
\end{equation}
$$\text{for} \;\; r\in (c_1t, \; (c_1+c_2)t] .$$

\vskip 0.2cm

$\bullet$ {\it Evaluation of} $\; \text{Pr} \left\{ \rho(t) < r, \; N_1(t)\ge 1, \; N_2(t)\ge 1 \right\}$. Since    
$$\left\{ N_1(t)\ge 1 \right\} = \left\{ \bold Z_1(t)\in \text{int} \; \bold B_{c_1t} \right\} = \left\{ R_1(t)\in (0, \; c_1t) \right\}, $$
$$\left\{ N_2(t)\ge 1 \right\} = \left\{ \bold Z_2(t)\in \text{int} \; \bold B_{c_2t} \right\} = \left\{ R_2(t)\in (0, \; c_2t)  \right\} ,$$
then, taking into account that $ R_1(t)$ and $ R_2(t)$ are independent, we have: 
\begin{equation}\label{main13}
\aligned 
& \text{Pr} \left\{ \rho(t) < r, \; N_1(t)\ge 1, \; N_2(t)\ge 1 \right\} \\
& = \text{Pr} \biggl\{ R_1^2(t) + R_2^2(t) - 2R_1(t) R_2(t) \cos(\theta(t)) < r^2, \; R_1(t)\in (0, \; c_1t) , \; R_2(t)\in (0, \; c_2t) \biggr\} \\
& = \text{Pr} \biggl\{ \eta(t) > \frac{R_1^2(t) + R_2(t)^2 - r^2}{2R_1(t) R_2(t)} , \; R_1(t)\in (0, \; c_1t) , \; R_2(t)\in (0, \; c_2t) \biggr\} \\
& = \int\limits_0^{c_1t} \int\limits_0^{c_2t} \text{Pr} \biggl\{ \eta(t) > \frac{\xi^2 + \zeta^2 - r^2}{2\xi\zeta} \; \biggl\vert \; R_1(t)=\xi, \;\; R_2(t)=\zeta \biggr\} \; \text{Pr} \bigl\{ R_1(t)\in d\xi \bigr\} \; \text{Pr} \bigl\{ R_2(t)\in d\zeta \bigr\} . 
\endaligned
\end{equation}
According to (\ref{llem3}), the conditional probability in the integrand is:
$$\aligned 
\text{Pr} & \biggl\{ \eta(t) > \frac{\xi^2 + \zeta^2 - r^2}{2\xi\zeta} \; \biggl\vert \;  R_1(t)=\xi, \; R_2(t)=\zeta  \biggr\} \\
& = \left\{ \aligned 1, \qquad\qquad & \text{if} \;\; \frac{\xi^2 + \zeta^2 - r^2}{2\xi\zeta} \in(-\infty, \; -1], \\
\frac{1}{\pi} \arccos{ \left( \frac{\xi^2 + \zeta^2 - r^2}{2\xi\zeta} \right)}, \qquad & \text{if} \;\; \frac{\xi^2 + \zeta^2 - r^2}{2\xi\zeta} \in (-1, \; 1], \\
0, \qquad\qquad & \text{if} \;\; \frac{\xi^2 + \zeta^2 - r^2}{2\xi\zeta} \in (1, \; +\infty). \endaligned \right. \\
\\
& = \left\{ \aligned 1, \qquad\qquad & \text{if} \;\; \xi + \zeta \le r , \\
\frac{1}{\pi} \arccos{ \left( \frac{\xi^2 + \zeta^2 - r^2}{2\xi\zeta} \right)}, \qquad & \text{if} \;\; \xi + \zeta > r \; \text{and} \; |\xi-\zeta|\le r , \\
0, \qquad\qquad & \text{if} \;\; |\xi - \zeta| > r . 
\endaligned \right.
\endaligned$$
Therefore, (\ref{main13}) becomes 
\begin{equation}\label{main14}
\aligned 
& \text{Pr} \left\{ \rho(t) < r, \; N_1(t)\ge 1, \; N_2(t)\ge 1 \right\} \\
& = \iint\limits_{\substack{\xi+\zeta\le r\\ 0<\xi<c_1t\\ 0<\zeta<c_2t}} \text{Pr} \bigl\{ R_1(t)\in d\xi \bigr\} \; \text{Pr} \bigl\{ R_2(t)\in d\zeta \bigr\} \\
& + \frac{1}{\pi} \iint\limits_{\substack{\xi+\zeta>r\\ |\xi-\zeta|\le r\\ 0<\xi<c_1t\\ 0<\zeta<c_2t}} \arccos{ \left( \frac{\xi^2 + \zeta^2 - r^2}{2\xi\zeta} \right)} \; \text{Pr} \bigl\{ R_1(t)\in d\xi \bigr\} \; \text{Pr} \bigl\{ R_2(t)\in d\zeta \bigr\} .
\endaligned
\end{equation}

For the sake of brevity we denote these integrals as follows:
\begin{equation}\label{mainINT}
\aligned 
\mathcal I_1(r,t) & = \iint\limits_{\substack{\xi+\zeta\le r\\ 0<\xi<c_1t\\ 0<\zeta<c_2t}} \text{Pr} \bigl\{ R_1(t)\in d\xi \bigr\} \; \text{Pr} \bigl\{ R_2(t)\in d\zeta \bigr\} , \\
\mathcal I_2(r,t) & = \iint\limits_{\substack{\xi+\zeta>r\\ |\xi-\zeta|\le r\\ 0<\xi<c_1t\\ 0<\zeta<c_2t}} \arccos{ \left( \frac{\xi^2 + \zeta^2 - r^2}{2\xi\zeta} \right)} \; \text{Pr} \bigl\{ R_1(t)\in d\xi \bigr\} \; \text{Pr} \bigl\{ R_2(t)\in d\zeta \bigr\} .
\endaligned
\end{equation}
Our goal is to evaluate integrals (\ref{mainINT}).

If $r\in(0, \; c_2t]$, then some simple geometric reasonings and formula (\ref{pprop6}) yield for the first integral in (\ref{mainINT}): 
$$\aligned 
\mathcal I_1(r,t) & = \frac{\lambda_1 \lambda_2 e^{-(\lambda_1+\lambda_2)t}}{c_1c_2} \int\limits_0^r \frac{\xi}{\sqrt{c_1^2t^2 - \xi^2}} 
\exp\left( \frac{\lambda_1}{c_1} \sqrt{c_1^2t^2 - \xi^2} \right) \\ 
& \hskip 4cm \times \biggl\{ \int\limits_0^{r-\xi} \frac{\zeta}{\sqrt{c_2^2t^2 - \zeta^2}} 
\exp\left( \frac{\lambda_2}{c_2} \sqrt{c_2^2t^2 - \zeta^2} \right) d\zeta \biggr\} d\xi . 
\endaligned$$
According to (\ref{intAB}), the interior integral in curl brackets is: 
$$\int\limits_0^{r-\xi} \frac{\zeta}{\sqrt{c_2^2t^2 - \zeta^2}} 
\exp\left( \frac{\lambda_2}{c_2} \sqrt{c_2^2t^2 - \zeta^2} \right) d\zeta = \frac{c_2}{\lambda_2} \biggl[ e^{\lambda_2t} - \exp\left( \frac{\lambda_2}{c_2} 
\sqrt{c_2^2t^2 - (r-\xi)^2} \right) \biggr] .$$
Then we obtain for the first integral in (\ref{mainINT}):
\begin{equation}\label{main15}
\aligned 
\mathcal I_1(r,t) & = \frac{\lambda_1}{c_1} \; e^{-\lambda_1t}  \int\limits_0^r \frac{\xi}{\sqrt{c_1^2t^2 - \xi^2}} 
\exp\left( \frac{\lambda_1}{c_1} \sqrt{c_1^2t^2 - \xi^2} \right) d\xi \\ 
& \quad - \frac{\lambda_1}{c_1} \; e^{-(\lambda_1+\lambda_2)t} \int\limits_0^r \frac{\xi}{\sqrt{c_1^2t^2 - \xi^2}} \exp\left( \frac{\lambda_1}{c_1} \sqrt{c_1^2t^2 - \xi^2} \right) \exp\left( \frac{\lambda_2}{c_2} \sqrt{c_2^2t^2 - (r-\xi)^2} \right) d\xi \\ 
& = 1 - \exp\left( -\lambda_1t + \frac{\lambda_1}{c_1} \sqrt{c_1^2t^2 - r^2} \right) \\ 
& \quad - \frac{\lambda_1}{c_1} \; e^{-(\lambda_1+\lambda_2)t} \int\limits_0^r \frac{\xi}{\sqrt{c_1^2t^2 - \xi^2}} \exp\left( \frac{\lambda_1}{c_1} \sqrt{c_1^2t^2 - \xi^2} \right) \exp\left( \frac{\lambda_2}{c_2} \sqrt{c_2^2t^2 - (r-\xi)^2} \right) d\xi ,
\endaligned
\end{equation}
$$\text{for} \;\; r\in(0, \; c_2t],$$ 
where in the last step we have used again formula (\ref{intAB}). 

In the same manner, if $r\in(c_2t, \; c_1t]$, then for the first integral in (\ref{mainINT}) we get: 
$$\aligned 
\mathcal I_1(r,t) & = \frac{\lambda_1 \lambda_2 e^{-(\lambda_1+\lambda_2)t}}{c_1c_2} \int\limits_0^{c_2t} \frac{\zeta}{\sqrt{c_2^2t^2 - \zeta^2}} 
\exp\left( \frac{\lambda_2}{c_2} \sqrt{c_2^2t^2 - \zeta^2} \right) \\ 
& \hskip 4cm \times \biggl\{ \int\limits_0^{r-\zeta} \frac{\xi}{\sqrt{c_1^2t^2 - \xi^2}} 
\exp\left( \frac{\lambda_1}{c_1} \sqrt{c_1^2t^2 - \xi^2} \right) d\xi \biggr\} d\zeta . 
\endaligned$$
Applying (\ref{intAB}), after some simple computations we arrive to the formula: 
\begin{equation}\label{main16}
\aligned 
\mathcal I_1(r,t) & = 1 - e^{-\lambda_2t} - \frac{\lambda_2}{c_2} e^{-(\lambda_1+\lambda_2)t} \\ 
& \qquad \times \int\limits_0^{c_2t} \frac{\zeta}{\sqrt{c_2^2t^2 - \zeta^2}} 
\exp\left( \frac{\lambda_2}{c_2} \sqrt{c_2^2t^2 - \zeta^2} \right) \; \exp\left( \frac{\lambda_1}{c_1} \sqrt{c_1^2t^2 - (r-\zeta)^2} \right) d\zeta ,
\endaligned
\end{equation}
$$\text{for} \;\; r\in(c_2t, \; c_1t] .$$

If $r\in(c_1t, \; (c_1+c_2)t]$, then for the first integral in (\ref{mainINT}) we have: 

$$\aligned 
\mathcal I_1(r,t) & = \frac{\lambda_1 \lambda_2}{c_1 c_2} \; e^{-(\lambda_1+\lambda_2)t}  \int\limits_0^{c_1t} \frac{\xi}{\sqrt{c_1^2t^2 - \xi^2}} 
\exp\left( \frac{\lambda_1}{c_1} \sqrt{c_1^2t^2 - \xi^2} \right) d\xi \\
& \hskip 4cm \times \int\limits_0^{c_2t} \frac{\zeta}{\sqrt{c_2^2t^2 - \zeta^2}} \exp\left( \frac{\lambda_2}{c_2} \sqrt{c_2^2t^2 - \zeta^2} \right) d\zeta \\
& - \frac{\lambda_1 \lambda_2}{c_1 c_2} \; e^{-(\lambda_1+\lambda_2)t}  \int\limits_{r-c_2t}^{c_1t} \frac{\xi}{\sqrt{c_1^2t^2 - \xi^2}} 
\exp\left( \frac{\lambda_1}{c_1} \sqrt{c_1^2t^2 - \xi^2} \right) \\
& \hskip 4cm \times \biggl\{ \int\limits_{r-\xi}^{c_2t} \frac{\zeta}{\sqrt{c_2^2t^2 - \zeta^2}} \exp\left( \frac{\lambda_2}{c_2} \sqrt{c_2^2t^2 - \zeta^2} \right) d\zeta \biggr\} d\xi. 
\endaligned$$
Applying again (\ref{intAB}), after some simple computations we finally obtain: 
\begin{equation}\label{main17}
\aligned 
& \mathcal I_1(r,t) \\ 
& = (1-e^{-\lambda_1t}) (1-e^{-\lambda_2t}) - e^{-(\lambda_1+\lambda_2)t} \left[ 1 - \exp\left( \frac{\lambda_1}{c_1} \sqrt{c_1^2t^2 - (r-c_2t)^2} \right) \right] \\
& - \frac{\lambda_1}{c_1} \; e^{-(\lambda_1+\lambda_2)t} \int\limits_{r-c_2t}^{c_1t} \frac{\xi}{\sqrt{c_1^2t^2 - \xi^2}} 
\exp\left( \frac{\lambda_1}{c_1} \sqrt{c_1^2t^2 - \xi^2} \right) \exp\left( \frac{\lambda_2}{c_2} \sqrt{c_2^2t^2 - (r-\xi)^2} \right) d\xi ,
\endaligned
\end{equation}
$$\text{for} \;\; r\in(c_1t, \; (c_1+c_2)t] .$$

Consider now the second integral in (\ref{mainINT}). If $r\in(0, \; c_2t]$ then the integration area represents either a 4-gon with the vertices $A=(r,\; 0), \; B=(0,\; r), \; C=(c_2t-r, \; c_2t), \; D=(c_2t+r, \; c_2t)$ in the case $c_2t+r\le c_1t$, or a 5-gon with the vertices $A=(r,\; 0), \; B=(0,\; r), \; C=(c_2t-r, \; c_2t), \; D=(c_1t, \; c_2t), \; E=(c_1t, \; c_1t-r)$ in the case $c_2t+r>c_1t$. In the first of these cases (the 4-gon) the second integral in (\ref{mainINT}) can be evaluated by the formula: 
\begin{equation}\label{main18}
\aligned
\mathcal I_2(r,t) & = \int\limits_0^r d\zeta \biggl\{ f_2^{ac}(\zeta,t) \int\limits_{r-\zeta}^{r+\zeta} \arccos{ \left( \frac{\xi^2 + \zeta^2 - r^2}{2\xi\zeta} \right)} \; f_1^{ac}(\xi,t) \; d\xi \biggr\} \\ 
& + \int\limits_r^{c_2t} d\zeta \biggl\{ f_2^{ac}(\zeta,t) \int\limits_{\zeta-r}^{\zeta+r} \arccos{ \left( \frac{\xi^2 + \zeta^2 - r^2}{2\xi\zeta} \right)} \; f_1^{ac}(\xi,t) \; d\xi \biggr\} , 
\endaligned
\end{equation}
$$\text{for} \;\; r\in(0, \; c_2t] \;\; \text{as} \;\; c_2t+r\le c_1t ,$$
where $f_1^{ac}(z,t), \; f_2^{ac}(z,t)$ are the absolutely continuous parts of the densities of processes $R_1(t), \; R_2(t)$, respectively, given by (\ref{pprop6}). 

In the second case (the 5-gon) the second integral in (\ref{mainINT}) can be evaluated as follows: 

\begin{equation}\label{main19}
\aligned
\mathcal I_2(r,t) & = \int\limits_0^r d\zeta \biggl\{ f_2^{ac}(\zeta,t) \int\limits_{r-\zeta}^{r+\zeta} \arccos{ \left( \frac{\xi^2 + \zeta^2 - r^2}{2\xi\zeta} \right)} \; f_1^{ac}(\xi,t) \; d\xi \biggr\} \\ 
& + \int\limits_r^{c_1t-r} d\zeta \biggl\{ f_2^{ac}(\zeta,t) \int\limits_{\zeta-r}^{\zeta+r} \arccos{ \left( \frac{\xi^2 + \zeta^2 - r^2}{2\xi\zeta} \right)} \; f_1^{ac}(\xi,t) \; d\xi \biggr\} \\
& + \int\limits_{c_1t-r}^{c_2t} d\zeta \biggl\{ f_2^{ac}(\zeta,t) \int\limits_{\zeta-r}^{c_1t} \arccos{ \left( \frac{\xi^2 + \zeta^2 - r^2}{2\xi\zeta} \right)} \; f_1^{ac}(\xi,t) \; d\xi \biggr\} , 
\endaligned
\end{equation}
$$\text{for} \;\; r\in(0, \; c_2t] \;\; \text{as} \;\; c_2t+r > c_1t .$$

Similarly, 
\begin{equation}\label{main20}
\aligned
\mathcal I_2(r,t) & = \int\limits_0^{c_2t} d\zeta \biggl\{ f_2^{ac}(\zeta,t) \int\limits_{r-\zeta}^{r+\zeta} \arccos{ \left( \frac{\xi^2 + \zeta^2 - r^2}{2\xi\zeta} \right)} \; f_1^{ac}(\xi,t) \; d\xi \biggr\} ,
\endaligned
\end{equation}
$$\text{for} \;\; r\in (c_2t, \; c_1t] \;\; \text{as} \;\; c_2t + r \le c_1t ,$$
and
\begin{equation}\label{main21}
\aligned
\mathcal I_2(r,t) & = \int\limits_0^{c_1t-r} d\zeta \biggl\{ f_2^{ac}(\zeta,t) \int\limits_{r-\zeta}^{r+\zeta} \arccos{ \left( \frac{\xi^2 + \zeta^2 - r^2}{2\xi\zeta} \right)} \; f_1^{ac}(\xi,t) \; d\xi \biggr\} \\
& + \int\limits_{c_1t-r}^{c_2t} d\zeta \biggl\{ f_2^{ac}(\zeta,t) \int\limits_{r-\zeta}^{c_1t} \arccos{ \left( \frac{\xi^2 + \zeta^2 - r^2}{2\xi\zeta} \right)} \; f_1^{ac}(\xi,t) \; d\xi \biggr\} ,
\endaligned
\end{equation}
$$\text{for} \;\; r\in (c_2t, \; c_1t] \;\; \text{as} \;\; c_2t + r > c_1t .$$

Finally, 
\begin{equation}\label{main22}
\aligned
\mathcal I_2(r,t) & = \int\limits_{r-c_2t}^{c_1t} d\xi \biggl\{ f_1^{ac}(\xi,t) \int\limits_{r-\xi}^{c_2t} \arccos{ \left( \frac{\xi^2 + \zeta^2 - r^2}{2\xi\zeta} \right)} \; f_2^{ac}(\zeta,t) \; d\zeta \biggr\} ,
\endaligned
\end{equation}
$$\text{for} \;\; r\in (c_1t, \;(c_1+c_2)t] .$$

Substituting expressions for the integrals $\mathcal I_1(r,t)$ and $\mathcal I_2(r,t)$ into (\ref{main14}) we obtain the fourth joint distribution in (\ref{main4}). 

\vskip 0.2cm

To combine all the formulas obtained above we should consider separately the alternatives A1, A2, A3, related to the possible cases $c_1>2c_2$, $c_1<2c_2$, $c_1=2c_2$, respectively.

\vskip 0.2cm 

{\it Alternative A1.} If $c_1>2c_2$ then, according to (\ref{mM}), $m(t)=c_2t, \; M(t)=(c_1-c_2)t$ and, therefore, the following partition of the interval $(0, \; (c_1+c_2)t]$ holds: 
$$(0, \; (c_1+c_2)t] = (0, \; c_2t] \cup (c_2t, \; (c_1-c_2)t] \cup ((c_1-c_2)t, \; c_1t] \cup (c_1t, \; (c_1+c_2)t] .$$
Under this partition, the function $G(r,t)$ represented by (\ref{G}) in the interval $(0, \; c_2t]$ emerges by putting into (\ref{main4}) the joint distributions given by formulas (\ref{main5}) (this term is zero in this case), (\ref{main9}) (with $\beta(r,t)= c_2t+r$ in this case), (\ref{mmain11}) (this term is also zero in this case), (\ref{main15}) and (\ref{main18}) (multiplied by $1/\pi$ in view of (\ref{main14})). 

The function $H_1(r,t)$ represented by (\ref{H1}) in the interval $(c_2t, \; (c_1-c_2)t]$ emerges by putting into (\ref{main4}) the joint distributions given by formulas (\ref{main5}) (this term is zero in this case), (\ref{main10}) (with $\beta(r,t)= c_2t+r$ in this case), (\ref{mmain11}) (this term is also zero in this case), (\ref{main16}) and (\ref{main20}) (multiplied by $1/\pi$ in view of (\ref{main14})). 

The function $Q(r,t)$ represented by (\ref{Q}) in the interval $((c_1-c_2)t, \; c_1t]$ emerges by putting into (\ref{main4}) the joint distributions given by formulas (\ref{main5}), (\ref{main10}) (with $\beta(r,t)= c_1t$ in this case), (\ref{main11}), (\ref{main16}) and (\ref{main21}) (multiplied by $1/\pi$ in view of (\ref{main14}). 

Finally, the function $U(r,t)$ represented by (\ref{U}) in the interval $(c_1t, \; (c_1+c_2)t]$ emerges by putting into (\ref{main4}) the joint distributions given by formulas (\ref{main5}), (\ref{main10}) (with $\beta(r,t)= c_1t$ in this case), (\ref{main12}), (\ref{main17}) and (\ref{main22}) (multiplied by $1/\pi$ in view of (\ref{main14})).

\vskip 0.2cm

{\it Alternative A2.} If $c_1<2c_2$, then, according to (\ref{mM}), $m(t)=(c_1-c_2)t, \; M(t)=c_2t$ and, therefore, we have the partition: 
$$(0, \; (c_1+c_2)t] = (0, \; (c_1-c_2)t] \cup ((c_1-c_2)t, \; c_2t] \cup (c_2t, \; c_1t] \cup (c_1t, \; (c_1+c_2)t] .$$
One can see that for such partition the functions $G(r,t), \; Q(r,t), \; U(r,t)$ are combined as above in alternative A1. The only difference is function $H_2(r,t)$ defined in the interval $((c_1-c_2)t, \; c_2t]$ by formula (\ref{H2}) which in this case emerges by putting into (\ref{main4}) the joint distributions (\ref{main5}), (\ref{main9}) (with $\beta(r,t)= c_1t$), (\ref{main11}), (\ref{main15}) and (\ref{main19}) (multiplied by $1/\pi$ in view of (\ref{main14})). 

\vskip 0.2cm

{\it Alternative A3.} If $c_1=2c_2$, then, according to (\ref{mM}), $m(t)=M(t)=c_2t$ and, therefore, we have the partition: 
$$(0, \; (c_1+c_2)t] = (0, \; c_2t] \cup (c_2t, \; c_1t] \cup (c_1t, \; (c_1+c_2)t] .$$
This means that in distribution function (\ref{Phi}) the interval $(m(t), \; M(t)]$, as well as function $H_k(r,t)$, vanish, while the remaining functions $G(r,t), \; Q(r,t), \; U(r,t)$ emerge as above. From this fact (\ref{mainA3}) follows. 

\vskip 0.2cm

The theorem is thus completely proved. $\square$ 

\bigskip 

We finish this section by presenting a result concerning the case of equal velocities. Suppose that both the random flights $\bold Z_1(t)$ and $\bold Z_2(t)$ are developing with the same speed $c_1=c_2=c$. In such case $0<\rho(t)<2ct$ with probability 1 for any $t>0$ and, therefore, the open interval $(0, \; 2ct)$ is the support of the distribution of process $\rho(t)$. 

\bigskip

{\bf Theorem 2.} {\it Under the condition $c_1=c_2=c$, the probability distribution function $\Phi(r,t)$ of the Euclidean distance $\rho(t)$ between 
two independent planar random flights $\bold Z_1(t)$ and $\bold Z_2(t)$ has the form:}
\begin{equation}\label{PhiEq}
\Phi(r,t) = \left\{ 
\aligned 0, \qquad & \text{if} \; r\in (-\infty, \; 0] ,\\
         V(r,t), \qquad & \text{if} \; r\in (0, \; ct] ,\\
         W(r,t), \qquad & \text{if} \; r\in (ct, \; 2ct] ,\\
         1, \qquad & \text{if} \; r\in (2ct, \; +\infty) ,
\endaligned \right. \qquad c_1 = c_2 = c, \;\; t>0,          
\end{equation}
{\it where functions $V(r,t), \; W(r,t),$ are given by the formulas:} 
\begin{equation}\label{V}
\aligned 
V(r,t) & = 1 - \exp\left( -\lambda_1t + \frac{\lambda_1}{c} \sqrt{c^2t^2 - r^2} \right) + \frac{e^{-(\lambda_1+\lambda_2)t}}{\pi} \arccos\left( 1 - \frac{r^2}{2c^2t^2} \right) \\
& + \frac{\lambda_1 e^{-(\lambda_1+\lambda_2)t}}{\pi c}  \int\limits_{ct-r}^{ct} \arccos{ \left( \frac{\xi^2 + c^2t^2 - r^2}{2ct \xi} \right)} \frac{\xi}{\sqrt{c^2t^2 - \xi^2}} \exp\left( \frac{\lambda_1}{c} \sqrt{c^2t^2 - \xi^2} \right) d\xi \\
& + \frac{\lambda_2 e^{-(\lambda_1+\lambda_2)t}}{\pi c}  \int\limits_{ct-r}^{ct} \arccos{ \left( \frac{\xi^2 + c^2t^2 - r^2}{2ct \xi} \right)} \frac{\xi}{\sqrt{c^2t^2 - \xi^2}} \exp\left( \frac{\lambda_2}{c} \sqrt{c^2t^2 - \xi^2} \right) d\xi \\
& - \frac{\lambda_1}{c} \; e^{-(\lambda_1+\lambda_2)t} \int\limits_0^r \frac{\xi}{\sqrt{c^2t^2 - \xi^2}} \exp\left( \frac{\lambda_1}{c} \sqrt{c^2t^2 - \xi^2} \right) \exp\left( \frac{\lambda_2}{c} \sqrt{c^2t^2 - (r-\xi)^2} \right) d\xi \\ 
& + \frac{1}{\pi} \int\limits_0^{ct} d\xi \biggl\{ f_1^{ac}(\xi,t) \int\limits_0^{ct} \arccos{ \left( \frac{\xi^2 + \zeta^2 - r^2}{2\xi\zeta} \right)} \; f_2^{ac}(\zeta,t) \; d\zeta \biggr\} \\
& - \frac{1}{\pi} \int\limits_r^{ct} d\xi \biggl\{ f_1^{ac}(\xi,t) \int\limits_0^{\xi-r} \arccos{ \left( \frac{\xi^2 + \zeta^2 - r^2}{2\xi\zeta} \right)} \; f_2^{ac}(\zeta,t) \; d\zeta  \biggr\} \\
& - \frac{1}{\pi} \int\limits_0^r d\xi \biggl\{ f_1^{ac}(\xi,t) \int\limits_0^{r-\xi} \arccos{ \left( \frac{\xi^2 + \zeta^2 - r^2}{2\xi\zeta} \right)} \; f_2^{ac}(\zeta,t) \; d\zeta  \biggr\} \\
& - \frac{1}{\pi} \int\limits_r^{ct} d\zeta \biggl\{ f_2^{ac}(\zeta,t) \int\limits_0^{\zeta-r} \arccos{ \left( \frac{\xi^2 + \zeta^2 - r^2}{2\xi\zeta} \right)} \; f_1^{ac}(\xi,t) \; d\xi \biggr\} , 
\endaligned
\end{equation}
\vskip 0.2cm
\begin{equation}\label{W}
\aligned
W(r,t) & = 1 - e^{-(\lambda_1+\lambda_2)t} \biggl[ \exp\left( \frac{\lambda_2}{c} \sqrt{c^2t^2 - (r-ct)^2} \right) - \frac{1}{\pi} \arccos\left( 1 - \frac{r^2}{2c^2t^2} \right) \biggr] \\
& + \frac{\lambda_1 e^{-(\lambda_1+\lambda_2)t}}{\pi c} \int\limits_{r-ct}^{ct} \arccos{ \left( \frac{\xi^2 + c^2t^2 - r^2}{2ct \xi} \right)} \; \frac{\xi}{\sqrt{c^2t^2 - \xi^2}} \; \exp\left( \frac{\lambda_1}{c} \sqrt{c^2t^2 - \xi^2} \right) \; d\xi \\
& + \frac{\lambda_2 e^{-(\lambda_1+\lambda_2)t}}{\pi c}  \int\limits_{r-ct}^{ct} \arccos{ \left( \frac{\xi^2 + c^2t^2 - r^2}{2ct \xi} \right)} \; \frac{\xi}{\sqrt{c^2t^2 - \xi^2}} \; \exp\left( \frac{\lambda_2}{c} \sqrt{c^2t^2 - \xi^2} \right) \; d\xi \\
& - \frac{\lambda_1}{c} \; e^{-(\lambda_1+\lambda_2)t} \int\limits_{r-ct}^{ct} \frac{\xi}{\sqrt{c^2t^2 - \xi^2}} 
\exp\left( \frac{\lambda_1}{c} \sqrt{c^2t^2 - \xi^2} \right) \exp\left( \frac{\lambda_2}{c} \sqrt{c^2t^2 - (r-\xi)^2} \right) d\xi \\
& + \frac{1}{\pi} \int\limits_{r-ct}^{ct} d\xi \biggl\{ f_1^{ac}(\xi,t) \int\limits_{r-\xi}^{ct} \arccos{ \left( \frac{\xi^2 + \zeta^2 - r^2}{2\xi\zeta} \right)} \; f_2^{ac}(\zeta,t) \; d\zeta \biggr\} , 
\endaligned
\end{equation}
{\it and functions $f_i^{ac}(z,t), \; i=1,2,$ are given by} (\ref{ContDens}) {\it with $c_1=c_2=c$}.  

\vskip 0.2cm

{\it Proof.} We give a sketch of the proof only. Similarly as above, one can show that the joint probabilities are:
\begin{equation}\label{PEq00}
\aligned 
\text{Pr} & \left\{ \rho(t) < r, \; N_1(t)=0, \; N_2(t)=0 \right\} \\
& \qquad = \left\{ \aligned 0, \qquad\qquad\qquad\qquad & \text{if} \;\; r\in (-\infty, \; 0],\\
\frac{e^{-(\lambda_1+\lambda_2)t}}{\pi} \arccos\left( 1 - \frac{r^2}{2c^2t^2} \right), \qquad & \text{if} \;\;
r\in (0, \; 2ct], \\
e^{-(\lambda_1+\lambda_2)t}, \qquad\qquad\qquad\qquad & \text{if} \;\; r\in (2ct, \; +\infty), \endaligned \right. 
\endaligned
\end{equation}
\vskip 0.2cm
\begin{equation}\label{PEq10LessCt}
\aligned 
\text{Pr} & \left\{ \rho(t) < r, \; N_1(t)\ge 1, \; N_2(t)=0 \right\} \\
& = \frac{\lambda_1 e^{-(\lambda_1+\lambda_2)t}}{\pi c}  \int\limits_{ct-r}^{ct} \arccos{ \left( \frac{\xi^2 + c^2t^2 - r^2}{2ct \xi} \right)} \; \frac{\xi}{\sqrt{c^2t^2 - \xi^2}} \; \exp\left( \frac{\lambda_1}{c} \sqrt{c^2t^2 - \xi^2} \right) \; d\xi , 
\endaligned
\end{equation}
$$\text{for} \;\; r\in (0, \; ct] ,$$
\vskip 0.2 cm
\begin{equation}\label{PEq10MoreCt}
\aligned 
& \text{Pr} \left\{ \rho(t) < r, \; N_1(t)\ge 1, \; N_2(t)=0 \right\} \\
& = e^{-\lambda_2t} \left[ 1 - \exp\left( -\lambda_1t+\frac{\lambda_1}{c} \sqrt{c^2t^2 - (r-ct)^2} \right) \right] \\
& + \frac{\lambda_1 e^{-(\lambda_1+\lambda_2)t}}{\pi c} \int\limits_{r-ct}^{ct} \arccos{ \left( \frac{\xi^2 + c^2t^2 - r^2}{2ct \xi} \right)} \; \frac{\xi}{\sqrt{c^2t^2 - \xi^2}} \; \exp\left( \frac{\lambda_1}{c} \sqrt{c^2t^2 - \xi^2} \right) \; d\xi ,
\endaligned
\end{equation}
$$\text{for} \;\; r\in (ct, \; 2ct] .$$
\vskip 0.2 cm
\begin{equation}\label{PEq01LessCt}
\aligned 
& \text{Pr} \left\{ \rho(t) < r, \; N_1(t)=0, \; N_2(t)\ge 1 \right\} \\
& = \frac{\lambda_2 e^{-(\lambda_1+\lambda_2)t}}{\pi c}  \int\limits_{ct-r}^{ct} \arccos{ \left( \frac{\xi^2 + c^2t^2 - r^2}{2ct \xi} \right)} \; \frac{\xi}{\sqrt{c^2t^2 - \xi^2}} \; \exp\left( \frac{\lambda_2}{c} \sqrt{c^2t^2 - \xi^2} \right) \; d\xi ,
\endaligned
\end{equation}
$$\text{for} \;\; r\in (0, \; ct] ,$$
\vskip 0.2 cm
\begin{equation}\label{PEq01MoreCt}
\aligned 
& \text{Pr} \left\{ \rho(t) < r, \; N_1(t)=0, \; N_2(t)\ge 1 \right\} \\
& = e^{-\lambda_1t} \left[ 1 - \exp\left( -\lambda_2t+\frac{\lambda_2}{c} \sqrt{c^2t^2 - (r-ct)^2} \right) \right] \\
& + \frac{\lambda_2 e^{-(\lambda_1+\lambda_2)t}}{\pi c}  \int\limits_{r-ct}^{ct} \arccos{ \left( \frac{\xi^2 + c^2t^2 - r^2}{2ct \xi} \right)} \; \frac{\xi}{\sqrt{c^2t^2 - \xi^2}} \; \exp\left( \frac{\lambda_2}{c} \sqrt{c^2t^2 - \xi^2} \right) \; d\xi ,
\endaligned
\end{equation}
$$\text{for} \;\; r\in (ct, \; 2t] ,$$
\vskip 0.2 cm
$$\text{Pr} \left\{ \rho(t) < r, \; N_1(t)\ge 1, \; N_2(t)\ge 1 \right\} = \mathcal I_1(r,t) + \frac{1}{\pi} \mathcal I_2(r,t),$$
where
\begin{equation}\label{INT1LessCt}
\aligned 
\mathcal I_1(r,t) & = 1 - \exp\left( -\lambda_1t + \frac{\lambda_1}{c} \sqrt{c^2t^2 - r^2} \right) \\ 
& \quad - \frac{\lambda_1}{c} \; e^{-(\lambda_1+\lambda_2)t} \int\limits_0^r \frac{\xi}{\sqrt{c^2t^2 - \xi^2}} \exp\left( \frac{\lambda_1}{c} \sqrt{c^2t^2 - \xi^2} \right) \exp\left( \frac{\lambda_2}{c} \sqrt{c^2t^2 - (r-\xi)^2} \right) d\xi ,
\endaligned
\end{equation}
$$\text{for} \;\; r\in(0, \; ct],$$ 
\vskip 0.2 cm
\begin{equation}\label{INT1MoreCt}
\aligned 
\mathcal I_1(r,t) & = (1-e^{-\lambda_1t}) (1-e^{-\lambda_2t}) - e^{-(\lambda_1+\lambda_2)t} \left[ 1 - \exp\left( \frac{\lambda_1}{c} \sqrt{c^2t^2 - (r-ct)^2} \right) \right] \\
& - \frac{\lambda_1}{c} \; e^{-(\lambda_1+\lambda_2)t} \int\limits_{r-ct}^{ct} \frac{\xi}{\sqrt{c^2t^2 - \xi^2}} 
\exp\left( \frac{\lambda_1}{c} \sqrt{c^2t^2 - \xi^2} \right) \exp\left( \frac{\lambda_2}{c} \sqrt{c^2t^2 - (r-\xi)^2} \right) d\xi ,
\endaligned
\end{equation}
$$\text{for} \;\; r\in(ct, \; 2ct] ,$$
\vskip 0.2 cm
\begin{equation}\label{INT2LessCt}
\aligned 
\mathcal I_2(r,t) & = \int\limits_0^{ct} d\xi \biggl\{ f_1^{ac}(\xi,t) \int\limits_0^{ct} \arccos{ \left( \frac{\xi^2 + \zeta^2 - r^2}{2\xi\zeta} \right)} \; f_2^{ac}(\zeta,t) \; d\zeta \biggr\} \\
& - \int\limits_r^{ct} d\xi \biggl\{ f_1^{ac}(\xi,t) \int\limits_0^{\xi-r} \arccos{ \left( \frac{\xi^2 + \zeta^2 - r^2}{2\xi\zeta} \right)} \; f_2^{ac}(\zeta,t) \; d\zeta  \biggr\} \\
& - \int\limits_0^r d\xi \biggl\{ f_1^{ac}(\xi,t) \int\limits_0^{r-\xi} \arccos{ \left( \frac{\xi^2 + \zeta^2 - r^2}{2\xi\zeta} \right)} \; f_2^{ac}(\zeta,t) \; d\zeta  \biggr\} \\
& - \int\limits_r^{ct} d\zeta \biggl\{ f_2^{ac}(\zeta,t) \int\limits_0^{\zeta-r} \arccos{ \left( \frac{\xi^2 + \zeta^2 - r^2}{2\xi\zeta} \right)} \; f_1^{ac}(\xi,t) \; d\xi \biggr\} , 
\endaligned
\end{equation}
$$\text{for} \;\; r\in(0, \; ct],$$ 
\vskip 0.2 cm
\begin{equation}\label{INT2MoreCt}
\aligned 
\mathcal I_2(r,t) & = \int\limits_{r-ct}^{ct} d\xi \biggl\{ f_1^{ac}(\xi,t) \int\limits_{r-\xi}^{ct} \arccos{ \left( \frac{\xi^2 + \zeta^2 - r^2}{2\xi\zeta} \right)} \; f_2^{ac}(\zeta,t) \; d\zeta \biggr\} ,
\endaligned
\end{equation}
$$\text{for} \;\; r\in (ct, \; 2ct].$$ 
In formulas (\ref{INT2LessCt}) and (\ref{INT2MoreCt}) functions $f_i^{ac}(z,t), \; i=1,2,$ are given by (\ref{ContDens}) with $c_1=c_2=c$. 

Now function (\ref{V}), defined in the interval $(0, \; ct]$, emerges by summing joint distributions given by formulas (\ref{PEq00}), (\ref{PEq10LessCt}), (\ref{PEq01LessCt}), (\ref{INT1LessCt}) and (\ref{INT2LessCt}) (the latter one should be multiplied by $1/\pi$). Function (\ref{W}), defined in the interval $(ct, \; 2ct]$, emerges by summing joint distributions given by formulas (\ref{PEq00}), (\ref{PEq10MoreCt}), (\ref{PEq01MoreCt}), (\ref{INT1MoreCt}) and (\ref{INT2MoreCt}) (the latter one should also be multiplied by $1/\pi$). The theorem is proved. $\square$ 

\bigskip

{\it Remark 1.} By means of tedious but simple computations one can check that for any $t>0$ and arbitrary speeds $c_1, c_2$ 
such that $c_1 > c_2$ the following limiting relations hold: 
\begin{equation}\label{Lim1}
\aligned 
\lim\limits_{r\to 0+0} G(r,t) & = 0, \\
\lim\limits_{r\to m(t)-0} G(r,t) & = G(m(t), t) = \lim\limits_{r\to m(t)+0} H_k(r,t) , \qquad k=1, 2, \\
\lim\limits_{r\to M(t)-0} H_k(r,t) & = H_k(M(t), t) = \lim\limits_{r\to M(t)+0} Q(r,t) , \qquad k=1, 2, \\
\lim\limits_{r\to c_1t-0} Q(r,t) & = Q(c_1t, t) = \lim\limits_{r\to c_1t+0} U(r,t), \\
\lim\limits_{r\to (c_1+c_2)t-0} U(r,t) & = U((c_1+c_2)t, t) = 1 .  
\endaligned
\end{equation}

Similarly, for any $t>0$ and arbitrary speed $c_1=c_2=c$ 
\begin{equation}\label{Lim2}
\aligned 
\lim\limits_{r\to 0+0} V(r,t) & = 0, \\
\lim\limits_{r\to ct-0} V(r,t) & = V(ct, t) = \lim\limits_{r\to ct+0} W(r,t) , \\
\lim\limits_{r\to 2ct-0} W(r,t) & = W(2ct, t) = 1 .
\endaligned
\end{equation}

Formulas (\ref{Lim1}) show that, for arbitrary speeds $c_1, c_2$ such that $c_1 > c_2$, the probability distribution function $\Phi(r,t)$ is continuous at the points $0,  \; m(t), \; M(t), \; c_1t, \; (c_1+c_2)t$ and, therefore, it is continuous in the whole interval (i.e. the support) $[0, \; (c_1+c_2)t]$. Similarly, relations (\ref{Lim2}) 
prove the continuity of $\Phi(r,t)$ in the interval $[0, \; 2ct]$ in the case of equal speeds $c_1=c_2=c$. This entirely accords with the structure of the distribution described above. 

\bigskip

{\it Remark 2.} The functions $G(r,t), H_k(r,t), Q(r,t), U(r,t), V(r,t), W(r,t)$ composing the probability distribution functions $\Phi(r,t)$ in formulas (\ref{Phi}), (\ref{mainA3}) and (\ref{PhiEq}), have fairly complicated forms and, obviously, cannot be computed explicitly. Therefore, these functions can be evaluated numerically only. One can see that each of them contains the terms of two kinds. The terms of first kind contain single integral and such terms are easily computable numerically (for given parameters $\lambda_i, c_i, \; i=1,2,$ and time parameter $t$) by means of standard package of mathematical programs (such as MATHEMATICA or MAPLE) and usual personal computer. The terms of second kind contain double integrals that cannot be evaluated directly. To overcome this difficulty, one may decompose interior integral into a series and then to take some finite number of its terms with their subsequent integration until the necessary accuracy is reached.

\enddocument